\RequirePackage{fix-cm}
\documentclass[fms,times,10pt]{amsart}
\usepackage[french,english]{babel}
\usepackage[T1]{fontenc}
\usepackage{latexsym,amssymb}
\usepackage{amsmath}
\usepackage{amsfonts}
\usepackage[section]{placeins}
\usepackage{url}
\usepackage{microtype}
\usepackage{booktabs, multirow}
\usepackage{amsthm, comment, mathrsfs, url, yfonts}
\usepackage{array}
\usepackage{color}  

\theoremstyle{plain}
\newtheorem{theorem}{Theorem}[section]
\newtheorem{proposition}[theorem]{Proposition}
\newtheorem{corollary}[theorem]{Corollary}
\newtheorem{lemma}[theorem]{Lemma}
\newtheorem{conjecture}[theorem]{Conjecture}

\theoremstyle{definition}

\newtheorem{definition}[theorem]{Definition}

\newcommand{\begpf}{\noindent{\bf Proof.}\enspace}
\newcommand{\epf}{{\ifhmode\unskip\nobreak\hfil\penalty50 \hskip1em
\else\nobreak\fi \nobreak\mbox{}\hfil\mbox{$\square$} \parfillskip=0pt
\finalhyphendemerits=0 \par\vskip5pt}}

\newcommand{\iso}{\smash{\mathop{\longrightarrow}\limits^{\sim}}}

\newcommand{\ol}{\overline}

\newcommand{\C}{\mathbf{C}}
\newcommand{\F}{\mathbf{F}}
\newcommand{\N}{\mathbf{N}}
\newcommand{\Q}{\mathbf{Q}}
\newcommand{\R}{\mathbf{R}}
\newcommand{\Z}{\mathbf{Z}}

\DeclareMathOperator{\GL}{GL}
\DeclareMathOperator{\PGL}{PGL}

\newcommand{\PSL}{\mathrm{PSL}}
\newcommand{\Ind}{\mathrm{Ind}}
\newcommand{\gr}{\mathrm{gr}}
\newcommand{\fil}{\mathrm{Fil}}
\newcommand{\ur}{\mathrm{ur}}

\newcommand{\ab}{\mathrm{ab}}
\newcommand{\Aut}{\mathrm{Aut}}
\renewcommand{\hom}{\mathrm{Hom}}

\newcommand{\Art}{\mathrm{Art}}
\newcommand{\Ext}{\mathrm{Ext}}

\newcommand{\ar}{\mathrm{ty}}
\newcommand{\mr}{\mathrm{gt}}
\newcommand{\nr}{\mathrm{un}}
\newcommand{\pr}{\mathrm{fl}}
\newcommand{\cg}{\mathrm{cg}}

\newcommand{\gal}{\mathrm{Gal}}
\newcommand{\tr}{\mathrm{tr}}
\newcommand{\ord}{\mathrm{ord}}
\renewcommand{\th}{\mathrm{th}}
\newcommand{\triv}{\mathrm{triv}}
\newcommand{\cyc}{\mathrm{cyc}}
\newcommand{\Nm}{\mathrm{Norm}}
\renewcommand{\max}{\mathrm{max}}
\newcommand{\AH}{\mathrm{AH}}
\newcommand{\dr}{\mathrm{dR}}

\newcommand{\CO}{\mathcal{O}}

\newcommand{\gn}{\mathfrak{n}}

\newcommand{\gp}{\mathfrak{p}}

\newcommand{\frob}{\mathsf{Frob}}

\newcommand{\Fpbar}{\ol{\F}_p}

\newcommand{\gen}{\alpha}
\newcommand{\bargen}{\alpha}

\newcommand{\Fr}{\mathrm{Fr}}
\newcommand{\cO}{\mathcal{O}}
\newcommand{\w}{\alpha}
\newcommand{\bN}{{\mathbf{N}}}

\title[Serre weights and wild ramification]{Serre weights and wild ramification in two-dimensional
Galois representations}

\author{\sc Lassina Demb\'el\'e}
\address{Mathematics Institute, University of Warwick, Coventry CV4 7AL, UK}
\email{lassina.dembele@gmail.com}
\author{\sc Fred Diamond}
\address{Department of Mathematics,
King's College London, London WC2R 2LS, UK}
\email{fred.diamond@kcl.ac.uk}
\author{\sc David P.\ Roberts}
\address{Division of Science and Mathematics, University of Minnesota Morris, Morris, MN 56267 USA}
\email{roberts@morris.umn.edu }


\begin{document}

\maketitle

\begin{abstract}
A generalization of Serre's Conjecture asserts that if $F$ is a totally real field, then certain characteristic $p$ representations
 of Galois groups over $F$ arise from Hilbert modular forms. Moreover it predicts the set of weights of such forms in terms 
 of the local behavior of the Galois representation at primes over~$p$. This characterization of the weights, which is formulated 
 using $p$-adic Hodge theory, is known under mild technical hypotheses if $p > 2$.  In this paper we give, under the 
 assumption that $p$ is unramified in $F$, a conjectural alternative description for the set of weights.
{Our approach is to use the Artin--Hasse exponential and local class field theory to construct bases for local Galois cohomology spaces in terms of which we identify subspaces that
should correspond to ones defined using $p$-adic Hodge theory. The resulting conjecture amounts
to an explicit description of wild ramification in reductions of certain crystalline Galois representations.
It enables the direct computation of the set of Serre weights of a Galois representation, which we illustrate with numerical examples.  A proof of this conjecture has been announced by 
Calegari, Emerton, Gee and Mavrides.} 
\end{abstract}

\section{Introduction}

A conjecture of Serre \cite{serre}, now a theorem of Khare and Wintenberger  %
 \cite{kw1, kw2},  asserts that if $p$ is prime and
$$\rho: G_\Q \to  \GL_2(\ol{\F}_p)$$
is a continuous, odd, irreducible representation, then $\rho$ arises
from a Hecke eigenform in the space $S_k(\Gamma_1(N))$
of cusp forms of some weight $k$ and level $N$.  
Serre in fact formulated a refined version of the conjecture specifying the minimal
such $k$ and $N$ subject to the constraints $k \ge 2$ and $p \nmid N$; a key point
is that the weight depends only on the restriction of $\rho$ to a decomposition group
at $p$, and the level on ramification away from $p$.
The equivalence between the weaker version of the conjecture and its refinement
was already known through the work of many authors for $p>2$, and finally settled
for $p=2$ as well by Khare and Wintenberger.

Buzzard, Jarvis and one of the authors \cite{bdj} considered
a generalization of Serre's conjecture to the setting of Hilbert modular forms
for a totally real number field $F$ and formulated an analogous refinement
for representations $\rho: G_F \to  \GL_2(\ol{\F}_p)$
assuming $p$ is unramified in $F$; versions without this assumption are
given in {\cite{schein:ijm}} %
and \cite{gee:type}.
The equivalence between the conjecture and its refinement was
proved, assuming $p>2$ and a Taylor--Wiles hypothesis on $\rho$,
in a series of papers by Gee and several sets of co-authors culminating
in \cite{gls} and \cite{gk}, with an alternative to the latter provided by Newton
\cite{newton}.  Generalizations to higher-dimensional Galois representations
have also been studied by Ash, Herzig and others beginning with \cite{as}; 
see \cite{ghs} for recent development.   

 One of the main difficulties in
even formulating refined versions of generalizations of Serre's conjecture
is in prescribing the weights; the approach taken in \cite{bdj} and subsequent
papers, at least if $\rho$ is wildly ramified at primes over $p$, is to do this
in terms of Hodge--Tate weights of crystalline lifts of $\rho$.
The main purpose of this paper is to make the recipe for the set of
weights more explicit.   In view of the connection between Serre
weights and crystalline lifts, this amounts to a conjecture in explicit
$p$-adic Hodge theory about wild ramification in reductions of
crystalline Galois representations.

We now explain this in more detail.  
Let $F$ be a totally real number field,
$\CO$ its ring of integers, $\gn$ a non-zero ideal of $\CO$, $S_F$
the set of embeddings $F\to \R$ and suppose
$\vec{k} \in \Z^{S_F}$ with all $k_\tau \ge 2$ and of the same
parity.  A construction completed by Taylor in \cite{rlt:inv} then
associates a $p$-adic  Galois representation to each Hecke eigenform
in the space of Hilbert modular cusp forms of weight $\vec{k}$ and level $\gn$.
One then expects that every continuous, irreducible, totally odd
$$\rho:G_F \to \GL_2(\ol{\F}_p)$$
is modular in the sense that it is arises as the reduction of such a
Galois representation.  One further expects that the prime-to-$p$ part of the
minimal level from which $\rho$ arises is its Artin conductor,
but the prediction of the possible weights is more subtle.
If $p$ is unramified in $F$, then a recipe is given in \cite{bdj}
in terms of the restrictions of $\rho$ to decomposition groups at
primes $\gp$ over $p$.   An interesting feature of this recipe not so
apparent over $\Q$ is the dependence of the conjectured weights
on the associated local extension class when the restriction at $\gp$
is reducible.  If
\begin{equation}\label{eqn:reducible}
\rho|_{G_{F_\gp}} \sim \left(\begin{array}{cc}\chi_1&*\\ 0 &
\chi_2\end{array}\right)\end{equation}
{for some characters $\chi_1,\chi_2: G_{F_{\gp}} \to \ol{\F}_p^\times$, then
the resulting short exact sequence
$$0 \to \ol{\F}_p(\chi_1)  \to V_\rho \to \ol{\F}_p(\chi_2)  \to 0$$
defines a class in 
$$c_\rho \in \Ext^1_{\Fpbar[G_K]}(\Fpbar(\chi_2),\Fpbar(\chi_1)) \cong 
 \Ext^1_{\Fpbar[G_K]}(\Fpbar,\Fpbar(\chi)) \cong
 H^1(G_K,\Fpbar(\chi))$$
where $K = F_\gp$ and $\chi = \chi_1\chi_2^{-1}$.
The class $c_\rho$ is well-defined up to a scalar in $\Fpbar^\times$, in the sense
that another choice of basis with respect to which $\rho|_{G_K}$ has the form
(\ref{eqn:reducible}) yields a non-zero scalar multiple of $c_\rho$.
Alternatively, we may view $c_\rho$ as the class in $H^1(G_K,\Fpbar(\chi))$
defined by the cocycle $z$ obtained by writing}
$${\chi_2^{-1} \otimes \rho|_{G_K} \sim  \left(\begin{array}{cc}\chi&z\\ 0 &
1\end{array}\right).}$$
The space $H^1(G_K,\Fpbar(\chi))$
has dimension at least $[K:\Q_p]$, with equality unless
$\chi$ is trivial or cyclotomic.   Whether $\rho$ is modular of a particular
weight depends on whether this extension class lies in a certain
distinguished subspace of  $H^1(G_K,\Fpbar(\chi))$ whose definition
relies on constructions from $p$-adic Hodge theory.  If $K \neq \Q_p$
and $\rho$ is wildly ramified at $\gp$, then the associated extension class
is a non-trivial element of a space of dimension at least two, making it
difficult to determine the set of weights without a more explicit description
of the distinguished subspaces.

We address the problem in this paper
by  using local class field theory and the Artin--Hasse exponential to give an
explicit basis for the space of extensions (Corollary~\ref{thm:basis}), in terms
of which we provide a conjectural alternate description of the distinguished
subspaces (Conjecture~\ref{conj:spaces}).   {We point out that a
related problem is considered by Abrashkin in \cite{abrashkin}; in particular,
the results of \cite{abrashkin} imply cases of our conjecture where the
distinguished subspaces can be described using the ramification
filtration on $G_K$.}

{An earlier version of this paper was posted on the arXiv in March 2016.  At the time, we reported that a proof of Conjecture~\ref{conj:spaces} under certain genericity hypotheses would be forthcoming in the Ph.D~thesis of Mavrides~\cite{mavrides}.   In fact, Conjecture~\ref{conj:spaces} has now been proved completely by Calegari, Emerton, Gee, and Mavrides in a preprint posted to the arXiv in August 2016 ~\cite{cegm}.  We remark that our restriction to the case where $K$ is unramified over $\Q_p$ is made essentially for simplicity.  The methods of this paper, and indeed of \cite{cegm}, are expected to apply to the general case where $K/\Q_p$ is allowed
to be ramified, but the resulting explicit description of the distinguished subspaces is likely to be much more complicated.} 

{The now-proved Conjecture~\ref{conj:spaces} immediately yields an alternate description of the set of Serre weights for $\rho$. Combining this with the predicted modularity of $\rho$ gives Conjecture~\ref{conj:serre2}, for which we have gathered extensive computational evidence.  Indeed the appeal of our description is that one can compute the set of Serre weights directly from $\rho$.  In this paper, we illustrate this computation systematically in several examples with $K/\Q_p$ quadratic and $p=3$.  A sequel paper~\cite{ddr} will support Conjecture~\ref{conj:serre2} via a much broader range of examples and  elaborate on computational methods. In particular, the examples provided in \cite{ddr} illustrate subtle features of the recipe for the weights arising only when $\chi$ is highly non-generic, with particular attention to the case $p=2$. Such examples were instrumental in leading us to Conjecture~\ref{conj:spaces} in its full generality.}

This paper is structured as follows:  In Section~\ref{subsec:weights} we
recall the general statement of the weight part of Serre's conjecture for
$F$ unramified at $p$.  In Sections~\ref{subsec:fil}, \ref{subsec:AH}
and~\ref{subsec:basis}, we study the space of extensions $H^1(G_K,\Fpbar(\chi))$
in detail, arriving at an explicit basis in terms of the Artin--Hasse exponential.
In Sections~\ref{subsec:dependence} and~\ref{subsec:subspaces}, we use
this basis to give our conjectural description of the distinguished
subspaces appearing in the definition of the set of Serre weights.
We illustrate this description in more detail in the quadratic case in
Section~\ref{sec:quadratic}, and with numerical examples for $p=3$
in Sections~\ref{sec:examples} and~\ref{sec:numerical}.  
We remark that aside from these examples
and the discussion of Serre's conjecture at the end of Sections~\ref{subsec:weights}
and~\ref{subsec:subspaces}, the setting for the paper is entirely local.    %

\section{Serre weights} \label{subsec:weights}
\subsection{Notation and general background}
Let $K$ be an unramified extension of $\Q_p$ with ring of integers $\CO_K$
and residue field $k$, and let $f = [K:\Q_p] = [k:\F_p]$.   We fix algebraic
closures $\overline{\Q}_p$ and $\overline{K}$ of $\Q_p$ and $K$,
and let $T$ denote the set of embeddings $K \to \ol{\Q}_p$.  We let
$\ol{\F}_p$ denote the algebraic closure $\F_p$ obtained as the residue field
of the ring of integers of $\ol{\Q}_p$, and we identify $T$ with
the set of embeddings $k \to \ol{\F}_p$ via the canonical bijection.

For a field $F$, we write $G_F$ for the absolute Galois group of $F$.
We let $I_K$ denote the inertia subgroup of $G_K$, i.e. the kernel of
the natural surjection $G_K \to G_k$.  We write $\frob$ for the absolute
(arithmetic) Frobenius elements on $k$ and on $\Fpbar$, and $\frob_K$ for the arithmetic Frobenius
element of $G_K/I_K \cong G_k$.   We let $\Art_K:K^\times \to G_K^\ab$ denote the
Artin map, normalized in the standard way, so
the image of any uniformizer of $K$ in $G_K/I_K$ is $\frob_K^{-1}$.

Recall that the fundamental character $\omega_f : G_K \to k^\times$ is 
defined by 
$$\omega_f(g) = g(\pi)/\pi \bmod \pi\CO_L$$
where $\pi$ is any root of $x^{p^f-1} + p = 0$ and $L = K(\pi) \subset \ol{K}$.
Then the composite of $\omega_f$ with the Artin map
$K^\times \to G_K^\ab \to k^\times$
is the homomorphism sending $p$ to $1$ and any element of
$\CO_K^\times$ to its reduction mod $p$.   Replacing $\pi$ by a root
of $x^{p^f-1} + up = 0$ for $u \in \CO_K^\times$ alters $\omega_f$
by an unramified character, so in fact $\omega_f|_{I_K}$ is independent
of the choice of uniformizer $up$ of $K$.
For each $\tau \in T$, we define the associated
fundamental character $\omega_\tau: I_K \to \ol{\F}_p^\times$
to be $\tau\circ\omega_f|_{I_K}$.

A {\em Serre weight} (for $\GL_2(K)$) is an irreducible
$\ol{\F}_p$-representation of $\GL_2(k)$.  Recall  that these are 
precisely the representations of the form
$$V_{\vec{d},\vec{b}} =
\bigotimes_{\tau\in
T}(\det{}^{d_\tau}\otimes_{k}\mathrm{Sym}^{b_\tau-
1}k^2)\otimes_{k,\tau}\ol{\F}_p,$$
where $d_\tau$, $b_\tau \in \Z$ and $1\le b_\tau\le p$ for each $\tau
\in T$.   Moreover we can assume that $0 \le d_\tau \le p-1$ for each
$\tau \in T$ and that $a_\tau < p - 1$ for some $p$,
in which case the resulting $(p^{f}-1)p^{f}$ representations
$V_{\vec{d},\vec{b}}$ are also inequivalent.

Let $\rho: G_K \to \GL_2(\Fpbar)$ be a continuous representation.
The next two subsections recall from \cite{bdj} the definition of the set $W(\rho)$
of Serre weights associated to $\rho$.

\subsection{Serre weights associated to a reducible representation $\rho$}
Suppose first that $\rho$ is reducible and write $\rho \sim
\begin{pmatrix}\chi_1&*\\ 0&\chi_2\end{pmatrix}$.  The isomorphism
class of $\rho$ is then determined by the ordered pair $(\chi_1,\chi_2)$
and a cohomology class
$c_\rho \in H^1(G_K,\Fpbar(\chi))$, where we set $\chi = \chi_1\chi_2^{-1}$.  We first define a set
\begin{equation} \label{eqn:red}
W'(\chi_1,\chi_2) = \left\{\,(V_{\vec{d},\vec{b}},J)\,\left|\, 
\begin{array}{cc}  J\subset T,&
\displaystyle \chi_1|_{I_K} = \prod_{\tau\in T}\omega_\tau^{d_\tau}\prod_{\tau\in
J}\omega_\tau^{b_\tau} \\&\\
\mbox{and}&
\displaystyle  \chi_2|_{I_K} = \prod_{\tau\in T}\omega_\tau^{d_\tau}\prod_{\tau\not\in
J}\omega_\tau^{b_\tau}\end{array}\right.\,\right\}.
\end{equation}

For each pair $(V,J) \in W'(\chi_1,\chi_2)$ we will
define a subspace $L_{V,J} \subset H^1(G_K,\Fpbar(\chi))$, %
 but we first need to recall the notion of labelled Hodge--Tate weights.

Recall that if $V$ is an $n$-dimensional vector space over $\ol{\Q}_p$
and $\rho: G_K \to \Aut_{\ol{\Q}_p}(V)$ is a crystalline (hence de Rham)
representation, then $D = D_{\dr}(V) = (B_{\dr} \otimes_{\Q_p} V)^{G_K}$
is a free module of rank $n$ over $K\otimes_{\Q_p} \ol{\Q}_p$
 endowed with an (exhaustive, separated)
decreasing filtration by (not necessarily free) $K\otimes_{\Q_p} \ol{\Q}_p$-submodules.
Writing $K \otimes_{\Q_p} \ol{\Q}_p \cong \prod_{\tau\in T} \ol{\Q}_p$, we have a corresponding
decomposition $D = \oplus_{\tau\in T} D_\tau$ where each $D_\tau$ is an $n$-dimensional
filtered vector space over $\ol{\Q}_p$.  For each $\tau\in T$,
the multiset of {\em $\tau$-labelled Hodge--Tate weights} of $V$
are the integers $m$ with multiplicity $\dim_{\ol{\Q}_p} \gr^{-m} D_\tau$.
In particular if $\psi: G_K \to \ol{\Q}_p^\times$ is a crystalline character, then
it has a unique $\tau$-labelled Hodge--Tate weight $m_\tau$ for each $\tau \in T$.
One finds that the vector $\vec{m} = (m_\tau)_{\tau \in T}$ determines $\psi$ up
to an unramified character, and that 
$\ol{\psi}|_{I_K} = \prod_{\tau\in T} \omega_\tau^{m_\tau}$.

Returning to the definition of $L_{V,J}$, suppose that $V = V_{\vec{d},\vec{b}}$.
Let $\tilde{\chi}_1$ be the crystalline lift of $\chi_1$
with $\tau$-labelled Hodge--Tate weight $d_\tau + b_\tau$ (resp.~$d_\tau$)
for $\tau \in J$ (resp.~$\tau\not\in J$) such that $\tilde{\chi}_1(\Art_K(p)) = 1$,
and similarly let $\tilde{\chi}_2$ be the crystalline lift of $\chi_2$
with $\tau$-labelled Hodge--Tate weight $d_\tau$ (resp.~$d_\tau+b_\tau$)
for $\tau \in J$ (resp.~$\tau\not\in J$) such that $\tilde{\chi}_2(\Art_K(p)) = 1$.
We then let $L_{V,J}'$ denote the set of extension classes associated
to reductions of crystalline extensions of $\tilde{\chi}_2$ by $\tilde{\chi}_1$.
We then set $L_{V,J} = L_{V,J}'$ except in the following two cases (continuing
to denote $\chi_1\chi_2^{-1}$ by $\chi$):
\begin{itemize}
\item If $\chi$ is cyclotomic, $\vec{b} = (p,\ldots,p)$ and $J=T$,
then $L_{V,J} = H^1(G_K,\Fpbar(\chi))$;
\item if $\chi$ is trivial and $J \neq T$, then 
$L_{V,J} = L_{V,J}' + H^1_\ur(G_K,\Fpbar(\chi))$
where $H^1_\ur(G_K,\Fpbar(\chi))$ is the set of unramified
homomorphisms $G_K \to \Fpbar$.
\end{itemize}

Finally we define $W(\rho)$ by the rule
\begin{equation} \label{eqn:Wprho}
V\in W(\rho) \quad\Longleftrightarrow\quad \mbox{\,$(V,J) \in
W'(\chi_1,\chi_2)$ and $c_\rho \in L_{V,J}$
for some $J \subset T$.}\end{equation}
Thus $V \in W(\rho)$ if and only if $c_\rho \in L_V$
where $L_V$ is defined as the union of the $L_{V,J}$ over 
$$S_V(\chi_1,\chi_2) = \{\,J\subset T\,|\, (V,J) \in W'(\chi_1,\chi_2)\,\}$$
(so $L_V$ depends on the ordered pair $(\chi_1,\chi_2)$,
and it is understood to be the empty set if $S_V(\chi_1,\chi_2) = \emptyset$).

We remark that $L_{V,J}$ has dimension
at least $|J|$, with equality holding unless $\chi = \chi_1\chi_2^{-1}$ is
trivial or cyclotomic {(see \cite[Lemma~3.12]{bdj})}.    Moreover in a typical situation (for example
if $\chi= \prod_{\tau \in T} \omega_\tau^{a_\tau}$ with
$1 < a_\tau < p-1$ for all $\tau$), the projection from $W'(\chi_1,\chi_2)$
to the set of subsets of $T$ is bijective, and the projection
to the set of Serre weights is injective (see \cite[Section~3.2]{bdj}).
In that case $W'(\chi_1,\chi_2)$ has
cardinality $2^f$ and hence so does that of $W(\rho)$ if $c_\rho = 0$.
On the other hand, one would expect that for ``most'' $\rho$, the class
$c_\rho$ does not lie in any of the proper subspaces $L_{V,J}$ for $J\neq T$,
so that $W(\rho)$ contains a single Serre weight.

It is not however true in general that the projection from $W'(\chi_1,\chi_2)$
to the set of Serre weights is injective, i.e., $S_V(\chi_1,\chi_2)$ may
have cardinality greater than $1$, in which case it is not immediate from
the definition of $L_V$ that it is a subspace of 
$H^1(G_K,\Fpbar(\chi))$.
However it is proved in \cite{gls} if $S_V(\chi_1,\chi_2) \neq \emptyset$
and $p > 2$, then there is an element $J_\max \in S_V(\chi_1,\chi_2)$
such that $L_V = L_{V,J_\max}$, so that $L_V$ is in fact a subspace.
{Indeed the proof of Theorem~9.1 of \cite{gls} shows that
if $p >2$ and $V = V_{\vec{d},\vec{b}}$, then $\rho$ has a crystalline lift 
with $\tau$-labelled Hodge--Tate weights $(\{d_\tau,d_\tau+b_\tau\})_{\tau\in T}$,
if and only if $S_V(\chi_1,\chi_2) \neq \emptyset$ and $c_\rho \in L_{V,J_\max}$.
It follows that $L_{V,J} \subset L_{V,J_\max}$ for all $J \in S_V(\chi_1,\chi_2)$
(using that $J_\max = T$ in the exceptional case where $\chi$ is cyclotomic and
$\vec{b} = (p,\ldots,p)$), and that $V \in W(\rho)$ if and only if}
$\rho$ has a crystalline lift with $\tau$-labelled
Hodge--Tate weights $(\{d_\tau,d_\tau+b_\tau\})_{\tau\in T}$.

The main aim of this paper is to use local class field theory
to give a more explicit description of $H^1(G_K,\Fpbar(\chi))$,
and to use this description to define subspaces which we conjecture
coincide with the $L_V$ (even for $p=2$).

\subsection{Serre weights associated to an irreducible representation $\rho$}
While the focus of this paper is on the case where $\rho$ is reducible,
for completeness  we recall the definition of $W(\rho)$ in the case where
$\rho$ is irreducible.   We let $K'$ denote the quadratic unramified extension
of $K$, $k'$ the residue field of $K'$, $T'$
the set of embeddings of $k'$ in $\ol{\F}_p$, and $\pi$
the natural projection $T' \to T$.  For $\tau' \in T'$, we let 
$\omega_{\tau'}$ denote the corresponding fundamental
character of $I_{K'} = I_K$.   Note that if $\rho$ is irreducible,
then it is necessarily tamely ramified and in fact induced from
a character of $G_{K'}$.   We define $W(\rho)$ by the rule:
\begin{equation}\label{eqn:irred}
\begin{array}{ccc} &&
\rho|_I\sim\prod_{\tau \in T}\omega_\tau^{d_\tau}
\begin{pmatrix}\prod_{\tau'\in J'}\omega_{\tau'}^{b_{\pi(\tau')}}&0\\
0&\prod_{\tau'\notin J'}\omega_{\tau'}^{b_{\pi(\tau')}}\end{pmatrix}\\
V_{\vec{d},\vec{b}}\in W(\rho) &\Longleftrightarrow& \\
&&\mbox{ for some $J'\subset T'$ such that $\pi:J'\iso T$.}\end{array}
\end{equation}
It is true in this case as well that $W(\rho)$ typically has
cardinality $2^f$ (see \cite[Section~3.1]{bdj}).  Moreover
the result of \cite{gls} characterizing $W(\rho)$ in terms of
reductions of crystalline representations (for $p > 2$)
holds in the irreducible case as well.

\subsection{The case $K=\Q_p$}  To indicate the level of complexity hidden in the general recipe
for weights,  we describe the set $W(\rho)$ more explicitly in the classical case $K=\Q_p$.  
Replacing $\rho$ by a twist, we can assume $\rho|_{I_{\Q_p}}$ has the form
$ \begin{pmatrix}\omega_2^a&{0}\\0&\omega_2^{p a}\end{pmatrix}$
or $ \begin{pmatrix}\omega^a&{*}\\0&1\end{pmatrix}$
for some $a$ with $1\le a \le p-1$ (where $\omega = \omega_1$ is the 
mod $p$ cyclotomic character). 
In the first case we find that $W(\rho) = \{V_{0,a},V_{a-1,p+1-a}\}$
(with the two weights coinciding if $a=1$).
In the second case we may further assume (after twisting) that
$\rho = \begin{pmatrix}\chi &{*}\\0& 1\end{pmatrix}$ for some character
$\chi: G_{\Q_p} \to \ol{\F}_p^\times$. Since the space $H^1(G_{\Q_p},\Fpbar(\chi))$
is one-dimensional unless $\chi$ is trivial or cyclotomic, one does not need
much information about the spaces $L_{V,J}$ in order to determine $W(\rho)$;
indeed all one needs is that:
\begin{itemize}
\item $L_{V,T}=H^1(G_{\Q_p},\Fpbar(\chi))$ unless $\chi$ is
cyclotomic and $V =V_{0,1}$;
\item $L_{V_{0,1},T}$ is the peu ramifi\'ee subspace if
$\chi$ is
cyclotomic, i.e., the subspace corresponding to
$\Z_p^\times \otimes \ol{\F}_p^\times$ under the Kummer isomorphism
$H^1(G_{\Q_p},\mu_p) \cong \Q_p^\times/(\Q_p^\times)^p$.
\item $L_{V,\emptyset}=0$ if $\chi \neq 1$.
\end{itemize}
It then follows (see \cite{bdj}) that
$$W(\rho) =  \left\{\begin{array}{ll}
\{V_{0,a}\},&
\mbox{if $1 < a < p- 1$ and $\rho$ is non-split,}\\
\{V_{0,a},V_{a,p-1-a}\},&
\mbox{if $1 < a < p - 2$ and $\rho$ is split,}\\
\{V_{0,p-2}, V_{p-2,p}, V_{p-2,1}\},&
\mbox{if $a = p - 2$, $p > 3$ and $\rho$ is split,}\\
\{V_{0,p-1}\},&
\mbox{if $a = p - 1$ and $p > 2$,}\\
\{V_{0,p}\},&
\mbox{if $a = 1$, $\chi=\omega$ and $\rho$
   is not peu ramifi\'ee,}\\
\{V_{0,p},V_{0,1},V_{1,p-2}\},&
\mbox{if $a = 1$, $p>3$ and $\rho$ is split,}\\
\{V_{0,3},V_{0,1},V_{1,3},V_{1,1}\},&
\mbox{if $a = 1$, $p=3$ and $\rho$ is split,}\\
\{V_{0,p},V_{0,1}\},&
\mbox{otherwise.}\end{array}\right.$$
We remark that the first case above is the most typical and the next one arises in
the setting of ``companion forms.''  The remaining cases take into account special
situations that arise when $\chi|_{I_K}$ or its inverse is trivial or cyclotomic.

\subsection{Serre's conjecture over totally real fields}
We now recall how Serre weights arise in the context of Galois representations
associated to automorphic forms.
Let $F$ be a totally real field in which $p$ is unramified.
Let $\CO_F$
denote its ring of integers and $S_p$ the set of primes of $\CO_F$
dividing $p$.
For each $\gp \in S_p$, we let $k_\gp = \CO_F/\gp$, $f_\gp = [k_\gp:\F_p]$
and $T_\gp$ the set of embeddings $\tau: k_\gp \to \ol{\F}_p$.
The irreducible $\ol{\F}_p$-representations of 
$\GL_2(\CO/p\CO) \cong \prod_{\gp\in S_\gp} \GL_2(k_\gp)$
are then of the form:
$V = \otimes_{\{\gp\in S_p\}} V_\gp$ where each $V_\gp$ is a Serre
weight for $\GL_2(F_\gp)$.

Suppose that $\rho:G_F \to \GL_2(\ol{\F}_p)$ is continuous, irreducible
and totally odd.   A notion of $\rho$ being {\em modular of weight} $V$ is
introduced in \cite{bdj}, where the following generalization
of Serre's Conjecture (from~\cite{serre}) is made:
\begin{conjecture} \label{conj:serre1}
{The representation} $\rho$ is modular of weight $V  = \otimes_{\{\gp\in S_p\}} V_\gp$
if and only if $V_\gp \in W(\rho|_{G_{F_\gp}})$ for all $\gp \in S_p$.
\end{conjecture}

We refer the reader to \cite{bdj} for the definition of modularity of weight
$V$ and its relation to the usual notion of weights of Hilbert modular forms.  
We just remark that $\rho$ is modular of {\em some} weight $V$
if and only if $\rho$ is {\em modular} in the usual sense that $\rho \cong \ol{\rho}_f$
for some Hilbert modular eigenform $f$, and that the set of weights for
which $\rho$ is modular determines the possible cohomological weights and local
behavior at primes over $p$ of those eigenforms (see \cite[Prop.~2.10]{bdj}).

Under the assumption that $\rho$ is modular (of some weight),
Conjecture~\ref{conj:serre1} can be viewed as the generalization
of the weight part of Serre's Conjecture and has been proved
under mild technical hypotheses (for $p>2$) in 
a series of papers by Gee and coauthors culminating in \cite{gls}, %
together with the results of either Gee and Kisin \cite{gk} or
Newton \cite{newton}.  Moreover their result holds without
the assumption that $p$ is unramified in $F$ using the description
of $W(\rho)$ in terms of reductions of crystalline representations.

Finally we remark that Conjecture~\ref{conj:serre1} is known in the
case $F = \Q$.  In this case the modularity of $\rho$ is a theorem of
Khare and Wintenberger~\cite{kw1,kw2}, and the weight part follows
from prior work of Gross, Edixhoven and others (see~\cite[Thm.~3.20]{bdj});
it amounts to the statement that if $2 \le k \le p +1$, then $\omega^d\rho$ arises
from a Hecke eigenform of weight $k$ and level prime to $p$ if
and only if $V_{d,k-1} \in W(\rho)$.

\section{The ramification filtration on cohomology}\label{subsec:fil}
In this section we use the upper numbering
of ramification groups to define filtrations on the Galois cohomology groups
parametrizing the extensions of characters under consideration.

\subsection{Definition of the filtration}
Continue to let $K$ denote a finite unramified extension of $\Q_p$ of degree $f$ with
residue field $k$, and
let $\chi:G_K \to \ol{\F}_p^\times$ be any character.  
Recall 
  {from \cite[IV.3]{serre_cl}} that $G_K$ has a decreasing filtration by closed subgroups $G_K^u$
where $G_K^{-1} = G_K$, $G_K^u = I_K$ for $-1 < u \le 0$, and 
$\bigcup_{u>0} G_K^u$ is the wild ramification
subgroup $P_K$.   We define an increasing filtration on $H^1(G_K,\Fpbar(\chi))$
by setting
$$\fil^s H^1(G_K,\Fpbar(\chi)) =
  \bigcap_{u > s-1} \ker (H^1(G_K,\Fpbar(\chi)) \to H^1(G_K^u,\Fpbar(\chi)))$$
for $s \in \R$.  Note that $\fil^s H^1(G_K,\Fpbar(\chi))= 0$ for $s < 0$, and that
$$ \fil^0 H^1(G_K,\Fpbar(\chi)) = \ker(H^1(G_K,\Fpbar(\chi)) \to H^1(I_K,\Fpbar(\chi))).$$
Let $z$ be a cocycle representing a class in $c \in H^1(G_K,\Fpbar(\chi))$.
Since $\chi|_{P_K}$ is trivial, the restriction of $z$ defines a homomorphism
$P_K \to \ol{\F}_p$;  so if $s \ge 1$, then $c \in \fil^s(H^1(G_K,\Fpbar(\chi)))$ if and only if 
$z(G_K^u) = 0$ for all $u > s - 1$.  In particular, $c \in \fil^1(H^1(G_K,\Fpbar(\chi)))$
if and only if $z(P_K) = 0$; since $H^1(I_K/P_K, \Fpbar(\chi)) = 0$, it follows that
$\fil^sH^1(G_K,\Fpbar(\chi)) = \fil^0H^1(G_K,\Fpbar(\chi))$ for $0 \le s \le 1$.

\subsection{Computation of the jumps in the filtration}
For any $s \in \R$, we set $\fil^{< s}(H^1(G_K,\Fpbar(\chi))) = \bigcup_{t < s}\fil^t(H^1(G_K,\Fpbar(\chi)))$.
Since $G_K^u = \bigcap_{v < u} G_K^v$, the compactness of $G_K$ and continuity
of the cocycle $z$ imply that in fact
$$\fil^{<s}(H^1(G_K,\Fpbar(\chi))) = \ker (H^1(G_K,\Fpbar(\chi)) \to H^1(G_K^{s-1},\Fpbar(\chi))).$$
We will now compute the jumps in the filtration, i.e., the dimension of
$$\gr^s(H^1(G_K,\Fpbar(\chi))) = \fil^s(H^1(G_K,\Fpbar(\chi)))/\fil^{<s}(H^1(G_K,\Fpbar(\chi)))$$
for every $s$ and $\chi$.

We must first introduce some notation.  Choose an embedding $\tau_0: k\to\ol{\F}_p$,
let $\tau_i = \tau_0\circ\frob^i$ where $\frob$ is the absolute Frobenius on $k$.
Recall that $\omega_f:G_K \to k^\times$ denotes the character defined by
$$\omega_f(g) = g(\pi)/\pi,$$
where $\pi$ is any root of $x^{p^f - 1} = - p$ in $\overline{K}$, and set 
$\omega_{f,i} = \omega_{\tau_i} =  \tau_i \circ\omega_f$ for $i=0,\ldots,f-1$.
We may then write $\chi|_{I_K} = \omega_{f,0}^n|_{I_K}$ where $n = \sum_{j=0}^{f-1} a_j p^j$
for integers $a_j$ satisfying $1 \le a_j \le p$ for $j = 0,\ldots, f-1$.  Moreover
this expression is unique if we further require (in the case that
$\chi|_{I_K}$ is the cyclotomic character) that some $a_j \neq p$ for some $j$.
We extend the definition of $a_j$ to all integers $j$ by setting
$a_j = a_{j'}$ if $j \equiv j' \bmod f$.   We define $(a_0,a_1,\ldots,a_{f-1})$ to be the 
{\em tame signature} of $\chi$;
thus the tame signature of $\chi$ is an element of the set
$$S =  \{\,1,2,\ldots,p\,\}^f - \{(p,p,\ldots,p)\}.$$
Define an action of $\gal(k/\F_p) = \langle \frob \rangle \cong \Z/f\Z$ on $S$
by the formula 
$$\frob \cdot (a_0,a_1,\ldots,a_{f-1}) = (a_{f-1},a_0,\ldots,a_{f-2}).$$
Note that if $\chi$ has tame signature $\vec{a}$, then
$\frob \circ \chi$ has tame signature $\frob(\vec{a})$, as does $\chi\circ \sigma$
where $\sigma$ is the (outer) automorphism of $G_K$ defined by conjugation
by a lift of $\frob \in \gal(k/\F_p) \cong \gal(K/\Q_p)$ to $G_{\Q_p}$.   We define
be the {\em period} of $\vec{a} \in S$ to be the
cardinality of its orbit under $\gal(k/\F_p)$, and the {\em absolute niveau} of $\chi$
to be the period of its tame signature.  (Note that the orbit of the tame signature of
$\chi$ under $\gal(k/\F_p)$ is independent of the choice of $\tau_0$.)

For $i = 0,\ldots,f-1$, we define
\begin{equation}\label{eqn:ni} n_i = \sum_{j=0}^{f-1} a_{i+j}p^j,\end{equation}
so that $n_0 \equiv n_ip^i\bmod (p^f-1)$ and $\chi|_{I_K} = \omega_{f,i}^{n_i}|_{I_K}$.

\begin{theorem}  \label{thm:slopes}
Let $d_s = \dim_{\ol{\F}_p} \gr^s(H^1(G_K,\Fpbar(\chi)))$ for $s \in \R$.
Then $d_s = 0$ unless $s =0$ or $1 < s \le  1+\frac{p}{p-1}$.  Moreover if
$d_s \neq 0$ and $1 < s < 1+ \frac{p}{p-1}$, then $s = 1 + \frac{m}{p^f-1}$ for
some integer $m$ not divisible by $p$.   More precisely, if $\chi$ has tame
signature $(a_0,a_1,\ldots,a_{f-1})$ of period $f'$
and the integers $n_i$ are defined by (\ref{eqn:ni}), then:
\begin{enumerate}
\item $d_0 = 1$ if $\chi$ is trivial and $d_0 = 0$ otherwise; 
\item if $1 < s < \frac{p}{p-1}$, then
$$d_s = \left\{\begin{array}{ll}
f/f', & \parbox{3in}{if $s = \frac{n_{i+k}}{p^f-1}$ for some $i,k$ such that $k > 0$, $a_i=p$,
$a_{i+1} = \cdots = a_{i+k-1} = p-1$ and $a_{i+k} \neq p-1$,}\\&\\
0,  & \mbox{otherwise;}\end{array}\right.$$
\item  if $\frac{p}{p-1} \le s < 1 + \frac{p}{p-1}$, then 
$$d_s = \left\{\begin{array}{ll}
f/f', & \mbox{if $s=1 + \frac{n_i}{p^f-1}$ for some $i$ such that $a_i \neq p$,}\\
0,  & \mbox{otherwise;}\end{array}\right.$$
\item $d_{1+\frac{p}{p-1}} = 1$ if $\chi$ is cyclotomic, and $d_{1+\frac{p}{p-1}} = 0$ otherwise.  
 \end{enumerate}
\end{theorem}
\begpf
We let $d_s'$ denote the value claimed for $d_s$ in the statement.
Note that if $1 < s < \frac{p}{p-1}$, then $d_s'$ is the number of $j \in R$ such that
$s=\frac{n_j}{p^f-1}$, where $R$ is the set of $j\in\{0,\ldots,f-1\}$ such that $a_j \neq p-1$ and
$(a_i,a_{i+1},\ldots,a_{j-1}) = (p,p-1,\ldots,p-1)$ for some $i$ with $j-f \le i < j$.
Moreover $R$ is in bijection with the set of $i \in \{0,\ldots,f-1\}$ such that $a_i = p$,
and if $j \in R$, then $1 < \frac{n_i}{p^f-1} < \frac{p}{p-1}$.  Therefore
$$\sum_{1 < s < \frac{p}{p-1}}\!\!\!\!\! d_s'   \,\, =\,\,  \#\, \{\, i \in \{0,\ldots,f-1\}\,|\, a_i = p \,\}.$$
Similarly if  $\frac{p}{p-1} \le s < 1 + \frac{p}{p-1}$, then $d_s'$ is the number
of $i \in \{0,\ldots,f-1\}$ such that $s=1 + \frac{n_i}{p^f-1}$ and $a_i \neq p$;
moreover if $a_i \neq p$, then $\frac{p}{p-1} \le 1 + \frac{n_i}{p^f-1} < 1 + \frac{p}{p-1}$, so
$$\sum_{\frac{p}{p-1} \le  s < 1 + \frac{p}{p-1}} \!\!\!\!\!\!\!\!\!\! d_s' \,\,  = \,\, \#\, \{\, i \in \{0,\ldots,f-1\}\,|\, a_i \neq p \,\}.$$
It follows that
$$\sum_{s\in\R} d_s' = \left\{\begin{array}{ll}
f+2,&\mbox{if $p=2$ and $\chi$ is trivial,}\\
f+1,&\mbox{if $p>2$ and $\chi$ is trivial or cyclotomic,}\\
f,&\mbox{otherwise.}\end{array}\right.$$
Therefore $\sum_{s\in \R} d_s' = \dim_{\ol{\F}_p} H^1(G_K,\Fpbar(\chi))
 = \sum_{s\in\R} d_s$, so it suffices to prove that $d_s' \le d_s$
 for all $s$, and we need only consider $s$ such that $d_s' > 0$.
 
 For $s=0$, the inflation-restriction exact sequence
 $$0 \to H^1(G_K/I_K,\Fpbar(\chi)^{I_K}) 
  \to H^1(G_K,\Fpbar(\chi)) \to H^1(I_K,\Fpbar(\chi))$$
 shows that $\gr^0 H^1(G_K,\Fpbar(\chi)) \cong H^1(G_K/I_K,\Fpbar(\chi)^{I_K})$
 has dimension $1$ if $\chi$ is trivial, and $0$ otherwise,
 so that $d_0 = d_0'$.   We may therefore assume that $s > 1$
 and that $m = (s-1)(p^f-1)$ is an integer.   Moreover
 either $0 < m < \frac{p(p^f-1)}{p-1}$ and $m$ is not divisible by $p$,
 or $m = \frac{p(p^f-1)}{p-1}$.
 
 Let $M = L(\pi)$ where $\pi^{p^f-1} = - p$ and $L$ is an unramified
 extension of $K$ of degree prime to $p$ such that $\chi|_{G_M}$
 is trivial;  thus $\chi = \mu\omega_{f,0}^{n_0}$ for
 some unramified character $\mu$ of $\gal(L/K)$.  Since $\gal(M/K)$
 has order prime to $p$, inflation-restriction gives
 $$H^1(G_K,\Fpbar(\chi)) \cong 
 H^1(G_M,\Fpbar(\chi))^{\gal(M/K)}
 = \hom_{\gal(M/K)}(G_M^\ab,\ol{\F}_p(\chi)),$$
 which we identify with 
 $$\hom_{\gal(M/K)}(M^\times,\ol{\F}_p(\chi))
  = \hom_{\gal(M/K)}(M^\times/(M^\times)^p,\ol{\F}_p(\chi))$$
  via the isomorphism $M^\times \cong W_M^\ab \subset G_M^\ab$
  of local class field theory.
 
Since $M$ is tamely ramified over $K$, we have $G_K^u \subset G_M$
for $u > 0$, and in fact
$G_K^u = G_M^{u(p^f-1)}$ by \cite[IV, Prop.15]{serre_cl}, which
maps onto $1 + \pi^{\lceil u(p^f-1) \rceil} \CO_M$ under the homomorphism
$W_M \to M^\times$ of local class field theory (see Cor.~3 to
Thm.~1 of \cite{serre_cl}).
Therefore a class in $H^1(G_K,\Fpbar(\chi))$ has trivial restriction
to $G_K^u$ for all $u > s-1$ (resp.~$G_K^{s-1}$)
if and only if the corresponding homomorphism
$M^\times/(M^\times)^p  \to \ol{\F}_p(\chi)$ factors through
$M^\times/(M^\times)^pU_{m+1}$ (resp.~$M^\times/(M^\times)^pU_m$),
where we write $U_t = 1+\pi^t\CO_M$ for a positive integer $t$.
It follows that
$$\gr^s(H^1(G_K,\Fpbar(\chi))) \cong \hom_{\gal(M/K)}(U_m/(U_m\cap (M^\times)^p)U_{m+1},
\ol{\F}_p(\chi)).$$

Now suppose that $m < \frac{p(p^f-1)}{p-1}$ and $m$ is not divisible by $p$.
Then we claim that $U_m \cap (M^\times)^p \subset U_{m+1}$.  Indeed suppose
that $v_\pi(x^p-1) = m$ for some $x \in M^\times$, and let $t = v_\pi(x-1)$.
Then $t > 0$ and writing $x = 1+ y\pi^t$ for some $y \in \CO_M^\times$, we have
$$x^p - 1 = (1+y\pi^t)^p - 1
      = py\pi^t  + \cdots + y^p \pi^{pt}.$$
So $m  \ge \min(t+p^f-1,tp)$, with equality unless $t+p^f-1 = tp$.
If $t + p^f - 1 > tp$, then $m = tp$ contradicts that $m$ is not divisible by $p$,
and if $t+p^f-1 \le tp$, then $t \ge \frac{p^f-1}{p-1}$ contradicts that $m < \frac{p(p^f-1)}{p-1}$.
This establishes the claim, from which it follows that 
$$\gr^s H^1(G_K,\Fpbar(\chi)) \cong \hom_{\gal(M/K)} (U_m/U_{m+1},\ol{\F}_p(\chi)).$$
Letting $l$ denote the residue field of $L$, the map $x \mapsto 1+x\pi^m$
induces a $\gal(M/K)$-equivariant isomorphism $l(\omega_f^m) \cong U_m/U_{m+1}$,
and the map $x \otimes 1 \mapsto (\sigma(x))_{\sigma}$ induces a
$\gal(M/K)$-equivariant isomorphism
$$l(\omega_f^m) \otimes_{\F_p} \ol{\F}_p \cong \bigoplus_{i=0}^{f-1}
  \left(\bigoplus_{\sigma\in S_i} \ol{\F}_p(\omega_{f,i}^m)\right)$$
where $S_i$ is the set of embeddings $l \to \ol{\F}_p$ restricting to $\tau_i$
and the action of $\gal(M/K)$ on  $\bigoplus_{\sigma\in S_i} \ol{\F}_p(\omega_{f,i}^m)$
is defined by $g((x_\sigma)_\sigma)= \omega_{f,i}^m(g)(x_{\sigma\circ g})_\sigma$.
Noting that $\bigoplus_{\sigma\in S_i} \ol{\F}_p \cong \Ind_{\gal(M/L)}^{\gal(M/K)}\ol{\F}_p$,
we see that
$$(U_m/U_{m+1}) \otimes_{\F_p} \ol{\F}_p \cong \bigoplus_{i=0}^{f-1} 
    \bigoplus_{\mu} \ol{\F}_p(\mu\omega_{f,i}^m),$$
where the second direct sum is over all characters $\mu:\gal(L/K) \to \ol{\F}_p^\times$.
Therefore $d_s$ is the number of $i$ such that $m \equiv n_i \bmod(p^f - 1)$.
The inequality $d_s' \le d_s$ is now immediate from the definition of $d_s'$.

Finally consider the case $s = 1 + \frac{p}{p-1}$, so $m = \frac{p(p^f-1)}{p-1}$;
we may assume $\chi$ is cyclotomic, and it suffices to prove that $d_s \ge 1$.
For $x \in U_{m+1}$, we see that $\exp(p^{-1}\log x)$ converges to a $p^\th$
root of $x$, so $U_{m+1} \subset (M^\times)^p$.   It follows that 
$\fil^s H^1(G_K,\Fpbar(\chi)) = H^1(G_K,\Fpbar(\chi))$.  Therefore it suffices to prove that
$\fil^{>s} H^1(G_K,\Fpbar(\chi)) \neq H^1(G_K,\Fpbar(\chi))$, i.e., that there is a class in
$H^1(G_K,\Fpbar(\chi))$ whose restriction to $G_K^{p/(p-1)}$ is non-trivial.
Since $G_K^{p/(p-1)} = G_{\Q_p}^{p/(p-1)}$, the diagram
$$\begin{array}{ccc}
H^1(G_{\Q_p},\Fpbar(\chi)) & \longrightarrow & H^1(G_{\Q_p}^{p/(p-1)},\Fpbar(\chi)) \\
&&\\
\downarrow &      & \parallel \\
&&\\
H^1(G_{K},\Fpbar(\chi)) & \longrightarrow & H^1(G_{K}^{p/(p-1)},\Fpbar(\chi)) \end{array}$$
reduces us to the case $K= \Q_p$, and we may further assume
$M = \Q_p(\pi) = \Q_p(\zeta_p)$.  We see in this case that if 
$x \in U_1$, then $x^p \in U_{p+1}$, so that $U_p \cap (M^\times)^p
\subset U_{p+1}$ (and in fact equality holds).   It follows that 
$$\gr^s(H^1(G_{\Q_p},\Fpbar(\chi)))
\cong \hom_{\gal(M/\Q_p)}(U_p/U_{p+1},\F_p(\chi)),$$
which is non-trivial (in fact one-dimensional) since $U_p/U_{p+1} \cong \F_p(\chi)$.
\epf

\subsection{Terminology associated with ramification}
Note that the dimensions $d_s$ in Theorem~\ref{thm:slopes}
are at most $1$ if $\chi$ has absolute niveau $f$,
in which case we say $\chi$ is {\em primitive}; otherwise we say $\chi$ is
{\em imprimitive}.  Thus $\chi$ is imprimitive if and only if
its tame signature $(a_0,a_1,\ldots,a_{f-1})$ has non-trivial rotational symmetry,
which is equivalent to $\chi$ extending to a character of $G_{K'}$ for some
proper subfield $K'$ of $K$ containing $\ol{\Q}_p$.

The statement of the theorem is also simpler if $a_i < p$ for all $i$, in which case
we say $\chi$ is {\em generic}; otherwise we say $\chi$ is {\em non-generic}.
Thus if $\chi$ is generic, then $d_s = 0$ if
$1 < s < \frac{p}{p-1}$ (by part (2) of the theorem); moreover $n_i \le p^f -1$
for all $i$, so we also have $d_s = 0$ if $2 < s < 1 + \frac{p}{p-1}$ (by part (3)
of the theorem).  To characterize the types of exceptional behavior arising in
extensions when $\chi$ is non-generic (or trivial or cyclotomic),
we introduce the following subspaces of $H^1(G_K,\Fpbar(\chi))$:
$$\begin{array}{rcl}
H^1_\nr(G_K,\Fpbar(\chi)) & = & \fil^0 H^1(G_K,\Fpbar(\chi)) = \fil^1 H^1(G_K,\Fpbar(\chi)); \\
H^1_\mr(G_K,\Fpbar(\chi)) & = & \fil^{< \frac{p}{p-1}} H^1(G_K,\Fpbar(\chi)); \\
H^1_\pr(G_K,\Fpbar(\chi)) & = & \fil^{\frac{p}{p-1}} H^1(G_K,\Fpbar(\chi)); \\
H^1_\cg(G_K,\Fpbar(\chi)) & = & \fil^2 H^1(G_K,\Fpbar(\chi)); \\
H^1_\ar(G_K,\Fpbar(\chi)) & = & \fil^{< 1 + \frac{p}{p-1}} H^1(G_K,\Fpbar(\chi)).\end{array}$$
We call $H^1_\nr(G_K,\Fpbar(\chi))$ the {\em unramified subspace} of $H^1(G_K,\Fpbar(\chi))$,
and we call $H^1_\mr(G_K,\Fpbar(\chi))$ (resp. $H^1_\pr(G_K,\Fpbar(\chi))$, $H^1_\cg(G_K,\Fpbar(\chi))$,
$H^1_\ar(G_K,\Fpbar(\chi))$)
the {\em gently} (resp.~{\em flatly, cogently, typically}) {\em ramified subspace}
of $H^1(G_K,\Fpbar(\chi))$.   We use the same terminology to describe the  cohomology
classes in these subspaces.

The following is immediate from Theorem~\ref{thm:slopes}:
\begin{corollary} \label{cor:slopes}
With the above notation, we have 
\begin{enumerate}
\item $H^1_\nr(G_K,\Fpbar(\chi)) = 0$ unless $\chi$
is trivial, in which case  $H^1_\nr(G_K,\Fpbar(\chi))$ has dimension $1$;
\item $H^1(G_K,\Fpbar(\chi))/H^1_\ar(G_K,\Fpbar(\chi)) = 0$ unless $\chi$ is cyclotomic,
in which case it has dimension $1$;
\item $H^1_\pr(G_K,\Fpbar(\chi))/H^1_\mr(G_K,\Fpbar(\chi))=0$ unless $\chi|_{I_K}$ is cyclotomic,
in which case it has dimension $f$;
\item $H^1_\ar(G_K,\Fpbar(\chi))/H^1_\nr(G_K,\Fpbar(\chi))$ has dimension $f$;
\item $H^1_\mr(G_K,\Fpbar(\chi))/H^1_\nr(G_K,\Fpbar(\chi))$
has dimension equal to the number of $i \in \{0,\ldots,f-1\}$ such that $a_i = p$;
\item $H^1_\cg(G_K,\Fpbar(\chi)) = H^1_\ar(G_K,\Fpbar(\chi))$ if $\chi$ is generic.
\end{enumerate}
\end{corollary}

Let $\rho: G_K \to \GL_2(\ol{\F}_p)$ be a reducible representation
of the form $\rho \sim \begin{pmatrix}\chi_1&*\\ 0&\chi_2\end{pmatrix}$
and $c_\rho$ an associated cohomology class.  For $s \ge 1$ (resp.~ $s > 1$),  we say
that $\rho$ has {\em slope at most $s$} (resp.~ {\em less than $s$}) if $G_K^u \subset \ker(\rho)$
for all $u > s - 1$ (resp.~ $u \ge s -1$), or equivalently if $c_\rho \in \fil^s(H^1(G_K,\Fpbar(\chi)))$
(resp.~ $\fil^{<s}H^1(G_K,\Fpbar(\chi))$).
Note that $\rho$ is (at most) tamely ramified if $c_\rho$ is unramified; we say that
$\rho$ is {\em gently} (resp.~{\em  flatly, cogently, typically}) {\em ramified} according to 
whether $c_\rho$ is.
We remark that \cite[2.1]{fontaine} shows that if $\rho$ arises from a
finite flat group scheme over $\CO_K$, then $c_\rho$ is flatly
ramified.  If $\chi$ is cyclotomic,
then our notion of flatly ramified coincides with Serre's notion 
of peu ramifi\'ee in \cite{serre} recalled above.

\section{The Artin--Hasse exponential}\label{subsec:AH}
In this section we establish some properties of the Artin--Hasse exponential
which strike us as having independent interest. %
{Recall from e.g. \cite[\S7.2]{robert} that the Artin--Hasse exponential is defined by a power series
with rational coefficients:
$$E_p(x) = \exp\left(\sum_{n\ge 0} \frac{x^{p^n}}{p^n}\right).$$
Here, as usual, $\exp(x) = \sum_{n\ge 0} (x^n/n!)$.}  Since $p$ is fixed
throughout, we will omit the subscript and simply write is as $E(x)$.
The denominators
of the coefficients of $E(x)$ are prime to $p$, so we
may regard $E(x) \in \Z_p[[x]]$, and hence as a function 
$E: B_r(0) \to B_r(1)$ for any $r < 1$, where $B_r(a)$ denotes
the open disk of radius $r$ of $a \in \C_p$.  

\subsection{First multiplicativity lemma for $E(x)$}
Let $l$ be a finite field and let $L$ be the field of fractions of $W(l)$.
For $a \in l$, let $[a] \in W(l)$ denote the Teichm\"uller lift of $a$.
The following lemma establishes the key property of $E(x)$ we will need;
we will use it to relate the additive structure of $l$ to the multiplicative structure
of tamely ramified extensions of $L$.

\begin{lemma}  \label{lem:AH} If $a,b \in l$ then
$E([a]x)E([b]x)E([a+b]x)^{-1} \in (W(l)[[x]]^\times)^p$.
\end{lemma}
\begpf  For $n\ge 0$, we define elements $\delta_n \in L$
inductively as follows:
$$\begin{array}{rcl}   \delta_0 &=& \frac{1}{p}\left([a]+[b]-[a+b])\right), \\
\delta_n &=& \frac{1}{p^n}\left(\varphi^n(\delta_0) - \sum_{i=0}^{n-1} p^i \delta_i^{p^{n-i}}\right)
  \qquad\mbox{for $n \ge 1$.}\end{array}$$
We claim  that $\delta_n \in W(l)$ for all $n \ge 0$.
The statement is clear for $n = 0$, so suppose that $n > 0$.
For $i = 0,\ldots,n-1$, we have $\delta_i^p \equiv \varphi(\delta_i) \bmod p$,
and therefore 
$$\delta_i^{p^{n-i}} = (\delta_i^p)^{p^{n-1-i}} \equiv \varphi(\delta_i)^{p^{n-1-i}} \bmod p^{n-i}.$$
By the definition of $\delta_{n-1}$ we have
$\varphi^{n-1}(\delta_0) = \sum_{i=0}^{n-1} p^i \delta_i^{p^{n-1-i}},$ so
$$\varphi^n(\delta_0) =  \sum_{i=0}^{n-1} p^i \varphi(\delta_i)^{p^{n-1-i}}
  \equiv \sum_{i=0}^{n-1} p^i \delta_i^{p^{n-i}} \bmod p^n,$$
 which gives the claim.
 
 Now consider the power series
 $$f(x) = \prod_{{i} \ge 0} E\left(\delta_i x^{p^i}\right),$$
 which converges in $W(l)[[x]]$ since $E\left(x^{p^i}\right) \equiv 1 \bmod x^{p^i}$.
 We claim that 
 $$E([a]x)E([b]x)E([a+b]x)^{-1} =  f(x)^p.$$
 We prove this working in $L[[x]]$, where
 $\exp(g(x))\exp(h(x)) = \exp(g(x)+h(x))$ for $g(x),h(x) \in xL[[x]]$, and
 therefore  $\exp(\sum_{i\ge 0} g_{i}(x) ) = \prod_{i\ge 0} \exp(g_i(x))$
 if $g_i(x) \in x^{p^i} L[[x]]$.
 Note first that we have
$$ E([a]x)E([b]x)E([a+b]x)^{-1}  = \exp\left(\sum_{n\ge 0} a_nx^{p^n}\right),$$
where $a_n = p^{-n}([a]^{p^n} + [b]^{p^n} - [a+b]^{p^n}) = p^{1-n}\varphi^n(\delta_0)$.
On the other hand 
$$\begin{array}{rcl} f(x)^p &=& \prod_{i \ge 0} E\left(\delta_i x^{p^i}\right)^p \\
 &=& \prod_{i\ge 0} \exp\left(p \sum_{m\ge 0} p^{-m}\delta_i^{p^m}x^{p^{i+m}}  \right)\\
 &=& \exp\left(\sum_{i,m\ge 0} p^{1-m}\delta_i^{p^m}x^{p^{i+m}}  \right)\\
&=&  \exp\left(\sum_{n\ge 0} b_nx^{p^n}\right),\end{array}$$
where
$$b_n = \sum_{i=0}^n p^{1+i-n}\delta_i^{p^{n-i}} = a_n.$$
This proves the claim and hence the lemma.
\epf

\subsection{Second multiplicativity lemma for $E(x)$}
We will also need the following property of $E(x)$, which will be used to
ensure that our constructions later are independent of various choices made.

\begin{lemma} \label{lem:delta}
If $\delta \in W(l)$, then
$$E(x)E((1+ p\delta)x)^{-1} \prod_{m\ge 0}  E(p\delta x^{p^m}) \in (W(l)[[x]]^\times)^p.$$
\end{lemma}
\begpf
We have
$$E((1+p\delta)x) E(x)^{-1} = \exp\left(\sum_{n\ge 0} a_nx^{p^n} \right)$$
where $\displaystyle{a_n = \sum_{i=1}^{p^n} \left(\begin{array}{c} p^n \\ i \end{array} \right) p^{i-n} \delta^i}$.
Note that $\left(\begin{array}{c} p^n \\ i \end{array} \right) p^{i-n}$ has valuation
$i - v_p(i)$ for $i=1,\ldots,p^n$, so that $a_n \equiv p\delta \bmod p^2$ if $p >2$,
and $a_n \equiv 2(\delta + \delta^2)\bmod 4$ if $p=2$.
On the other hand
$$\prod_{m\ge 0}  E(p\delta x^{p^m})   = \exp\left(\sum_{n\ge 0} b_nx^{p^n} \right)$$
where
$\displaystyle{b_n = \sum_{j=0}^n p^{p^j - j} \delta^{p^j}}$.  Setting  $c_n = p^{-1}(b_n - a_n)$
gives $\exp(c_n x^{p^n}) \in W(l)[[x]]^\times$ since $c_n \in pW(l)$, and  
$$E(x)E((1+ p\delta)x)^{-1} \prod_{m\ge 0}  E(p\delta x^{p^m})
    =  (\prod_{n\ge 0} \exp(c_n x^{p^n}))^p.$$
\epf

\subsection{Homomorphisms induced by $E(x)$}
Suppose now that $M$ is a subfield of $\C_p$ containing $L$,
and $\alpha \in M$ is such that $|\alpha| < 1$.  Note that
Lemma~\ref{lem:AH} yields a homomorphism
$\varepsilon: l \to \CO_L[[x]]^\times \otimes \F_p$ defined by
$\varepsilon([a]) = E([a]x) \otimes 1$.  
We can therefore define a homomorphism
\begin{equation} \label{def:vareps}
\varepsilon_{\alpha}: l\otimes \Fpbar  \to \CO_M^\times \otimes \Fpbar
\end{equation}
as the extension of scalars of the composite of $\varepsilon$ with the 
{multiplicative homomorphism $\CO_L[[x]]^\times \to \CO_M^\times$}
 induced by evaluation at $\alpha$, so that
$\varepsilon_{\alpha}(a\otimes b) = E([a]\alpha) \otimes b$.

In addition to properties of $\varepsilon_\alpha$ derived from Lemmas~\ref{lem:AH}
and~\ref{lem:delta}, we will also need the following:
\begin{lemma} \label{lem:shift}   If $|\alpha| < p^{-1/p(p-1)}$, then
$$\varepsilon_{\alpha^p} \circ \frob = \varepsilon_{-p\alpha}$$
where $\frob$ is the automorphism of $l\otimes \Fpbar$ induced by the absolute
Frobenius on $l$.
\end{lemma}
\begpf It suffices to prove that if $\beta \in \CO_M$ is such that
$|\beta| {<} p^{-1/p(p-1)}$,
then $E(\beta^p) E(-p\beta)^{-1} \in (\CO_M^\times)^p$.  On the one hand
we have $E(\beta^p) = \exp(-p\beta)E(\beta)^p$.  On the other hand, setting
$\gamma = \sum_{n\ge 1} p^{p^n - n - 1} (-\beta)^{p^n}$, we see that $|\gamma|
 < p^{-1/(p-1)}$, so $\exp(\gamma)$ converges to an element of $\CO_M^\times$ such that
 $$E(-p\beta) = \exp(-p\beta + p\gamma) = \exp(-p\beta)\exp(\gamma)^p.$$
 \epf

\section{A basis for the cohomology}  \label{subsec:basis}
We return to the setup of Section~\ref{subsec:fil}, so $K$ is an unramified
extension of $\Q_p$ of degree $f$ with residue field $k$, $T$ is the set
of embeddings $k \to \Fpbar$, 
and $\chi$ is a character $G_K \to \Fpbar^\times$.
We will use a homomorphism of the form (\ref{def:vareps})
 to construct an explicit basis
for $H^1(G_K,\Fpbar(\chi))$.

Let $M$ be a tamely ramified abelian extension of $K$ such that $\chi|_{G_M}$ is trivial.
We assume $M$ is of the form $L(\pi)$ where $L$ is an unramified extension of
$K$ of degree prime to $p$ and $\pi$ is a uniformizer 
of $M$ such that  $\pi^e \in K^\times$ where the total ramification degree
$e$ of $M$ divides $p^f-1$.  We thus allow $M$ and $\pi$ to have a more
general form than in the proof of Theorem~\ref{thm:slopes}, but note that we
still have 
$$H^1(G_K,\Fpbar(\chi)) \cong \hom_{\gal(M/K)}(M^\times, \Fpbar(\chi)),$$
which we identify with the $\Fpbar$-dual of the vector space
$$ U_\chi = \left( M^\times \otimes \Fpbar(\chi^{-1}) \right)^{\gal(M/K)}.$$
Our explicit basis for $H^1(G_K,\Fpbar(\chi))$ will be defined as the dual basis
to one we construct for $U_\chi$.

\subsection{Definition of $u_i$}
As in Section~\ref{subsec:fil} we choose an embedding
$\tau_0: k \to \Fpbar$ and let $\tau_i = \tau_0\circ\frob^i$, let
$(a_0,\ldots,a_{f-1})$ be the tame signature of $\chi$ and define
the integers $n_i$ by (\ref{eqn:ni}), so that
$\chi|_{I_K} = \omega_{f,i}^{n_i}|_{I_K}$.
We define $\tau_i$, $a_i$ and $n_i$ for all $i \in \Z$ by requiring
that they depend only on $i\bmod f$.

Since $\chi|_{I_K}$ has
order dividing $e$, we see that $n_i$ is divisible by $(p^f-1)/e$
for all $i$. 
Letting  $\omega_\pi : \gal(M/K) \to \mu_e(K) \subset K^\times$ be the
character defined by $\omega_\pi(g) = g(\pi)/\pi$, we see that
$\overline{\omega}_\pi|_{I_K} = \omega_f|_{I_K}^{(p^f-1)/e}$, so that
$$\chi|_{I_K}  = (\tau_i\circ \overline{\omega}_\pi)|_{I_K}^{en_i/(p^f-1)}.$$

We now define an embedding $\tau_i'$ and
an integer $n_i'$ for each $i$.   If $a_{i+1} \neq p$, then we set $\tau_i' = \tau_{i+1}$
and $n'_i = en_{i+1}/(p^f-1)$.  If $a_{i+1} = p$, then we let $j$ be the least
integer greater than $i$ such that $a_{j+1} \neq p - 1$; thus 
$(a_{i+1},a_{i+2},\ldots,a_j) = (p,p-1,\ldots,p-1)$, but  $a_{j+1} \neq p - 1$.
We then set $\tau_i' = \tau_{j+1}$ and $n_i' = en_{j+1}/(p^f-1) - e$.
Note that for each $i$ we have $n_i' > 0$ and
$$\chi = \mu (\tau'_i \circ \overline{\omega}_\pi)^{n_i'}$$
for some unramified character $\mu: \gal(L/K) \to \Fpbar^\times$
independent of $i$.

Recall that we have an isomorphism
$$l\otimes \Fpbar \cong  \bigoplus_{\tau \in T}( l \otimes_{k,\tau}  \Fpbar)$$
defined by the natural projection on to each component.
By the Normal Basis Theorem, $l$ is free of rank one over $k[\gal(L/K)]
 = k[\gal(l/k)]$, so each component in the above decomposition
 is free of rank one over $\Fpbar[\gal(L/K)]$.  It follows that for
 each embedding $\tau \in T$, the $\mu$-eigenspace
 $$  \Lambda_{\tau,\mu}  =
 \{\, a \in l \otimes_{k,\tau}  \Fpbar \, | \, \mbox{$ga = (1\otimes \mu(g))a$ for all $g\in \gal(L/K)$}\,\}$$
 is one-dimensional over $\Fpbar$.
Let $\lambda_{\tau,\mu}$ be a non-zero element of $\Lambda_{\tau,\mu}$.
We now define
\begin{equation}
\label{eqn:units}
u_i = \varepsilon_{\pi^{n_i'}} (\lambda_{\tau_i',\mu})  \in \CO_M^\times \otimes \Fpbar
\end{equation}
for $i = 0,\ldots,f-1$, where $\varepsilon_{\pi^{n_i'}}$ is defined by (\ref{def:vareps}).

For $g \in \gal(M/K)$, $a \in l$, $n\in \Z$, we have 
$$g (E([a]\pi^n)) = E(g([a]\pi^n)) = E(\omega_\pi(g)^n g([a]) \pi^n),$$
so that
$$g(\varepsilon_{\pi^n}(\lambda))
      = \varepsilon_{\pi^n}((\overline{\omega}_\pi(g)^n \otimes 1) g(\lambda))$$
for all $\lambda \in l \otimes \Fpbar$.   Since
$$(\overline{\omega}_\pi(g)^{n_i'} \otimes 1) g(\lambda_{\tau_i',\mu})
  = (\overline{\omega}_\pi(g)^{n_i'} \otimes \mu(g)) \lambda_{\tau_i',\mu}
   = (1 \otimes \chi(g)) \lambda_{\tau_i',\mu},$$
we conclude that $gu_i = (1\otimes\chi(g))u_i$ for all $g \in \gal(M/K)$.
We can therefore view $u_i$ as an element of $U_\chi$.

\subsection{Definition of $u_{\rm triv}$ and $u_{\rm cyc}$}
We now define additional classes in $U_\chi$ in the case that $\chi$
is trivial or cyclotomic.  Note that if $g\in \gal(M/K)$, then
$g(\pi) = \omega_\pi(g)\pi$ and $\omega_\pi(g) \in \mu_e(K)$,
which is contained in $(M^\times)^p$.  It follows that
$g(\pi \otimes 1) = \pi \otimes 1$ in $M^\times \otimes \Fpbar$,
so that 
\begin{equation}\label{eqn:triv}
 u_\triv = \pi \otimes 1
\end{equation}
defines an element of $U_\chi$ for
the trivial character $\chi$.

If $\chi$ is cyclotomic then the assumption that $\chi|_{G_M}$ is trivial
ensures that $\Q_p(\zeta_p) \subset M$; in particular $e$ is divisible by $p-1$.
We now determine which elements $\alpha \in 1 + \pi^{ep/(p-1)}\CO_M$ are
$p^\th$-powers.  Recall that $\Q_p(\zeta_p) = \Q_p(\delta)$ where 
$\delta^{p-1} = -p$, so we may write $\alpha$
in the form $1 + \beta\delta^p$ with $\beta \in \CO_M$.  
We claim that $\alpha = 1+\beta\delta^p \in (M^\times)^p$ if and only
if $\tr_{l/\F_p} \overline{\beta} = 0$.     Suppose first that
$\tr_{l/\F_p}\overline{\beta} = 0$. We can then write
$\overline{\beta} = \overline{\gamma} - \overline{\gamma}^p$ for
some $\overline{\gamma} \in l$.   {(This follows for example from the fact that
$H^1(\gal(l/\F_p),l) = 0$ since $l$ is free over $\F_p[\gal(l/\F_p)]$ by the
Normal Basis Theorem.  Alternatively, note that $\tr_{l/\F_p}$ is surjective since
$l$ is separable over $\F_p$, so counting dimensions shows that
$l \to l \to \F_p$ is exact, where the maps are $1-\frob$ and $\tr_{l/\F_p}$.)
We can therefore write 
$\beta \equiv \gamma - \gamma^p \bmod \pi\CO_M$ for some
$\gamma \in \CO_M$, so that
$$ \alpha(1 + \gamma\delta)^p \equiv 1 + (\beta - \gamma + \gamma^p)\delta^p
      \equiv 1 \bmod \pi^{1 + \frac{ep}{p-1}}\CO_M.$$
Since $\exp(p^{-1}\log x)$ converges
to a $p^\mathrm{th}$ root of $x$ if $|x-1| < p^{-p/(p-1)}$,}
 it follows that
$\alpha(1+\gamma\delta)^p \in  (\CO_M^\times)^p$,
and hence that $\alpha \in (\CO_M^\times)^p$.
Suppose conversely that $\alpha = 1 + \beta\delta^p \in (M^\times)^p$.
Then considering valuations as in the proof of Theorem~\ref{thm:slopes}, 
we see that $\alpha = (1+\gamma\delta)^p$
for some $\gamma \in \CO_M$.  Since 
$$(1+\gamma\delta)^p  \equiv  1 + (\gamma^p - \gamma)\delta^p \bmod \delta^{p+1}\CO_M,$$
we deduce that $\beta \equiv \gamma^p - \gamma \bmod \pi$, and hence
that $\tr_{l/\F_p}(\overline{\beta}) = 0$.  This completes the proof of the claim.
Now choose any $b \in l$ such that $\tr_{l/\F_p} b \neq 0$, and define
\begin{equation}
\label{eqn:cyc}
u_\cyc = \varepsilon_{\delta^p}(b\otimes 1) = E([b]\delta^p) \otimes 1,
\end{equation}
which is a non-trivial element of $\CO_M^\times \otimes \Fpbar$ by the above claim.
Moreover since  $\tr_{l/\F_p}(gb) = \tr_{l/\F_p}(b)$ for all
 $g \in \gal(M/K)$, it also follows from the claim and Lemma~\ref{lem:AH}
 that $E(g([b])\delta^p)E([b]\delta^p)^{-1} \in (\CO_M^\times)^p$,
so $\varepsilon_{\delta^p}(gb\otimes 1) = 
\varepsilon_{\delta^p}(b\otimes 1)$.   Since $g(\delta) = [\chi(g)]\delta$,
we see as in the construction of the elements $u_i$ that 
$$g(u_\cyc) = 
\varepsilon_{\delta^p}(\chi(g)gb \otimes 1)
  = (1\otimes\chi(g))\varepsilon_{\delta^p}(gb \otimes 1)
   =  (1\otimes\chi(g))u_\cyc,$$
and therefore that $u_\cyc \in U_\chi$.

\subsection{Bases for $U_\chi$ and {$H^1(G_K,\Fpbar(\chi))$} }
\begin{theorem}  \label{prop:basis}
Let $B$ denote the subset of $U_\chi$ consisting of the
elements $u_i$ for $i = 0,\ldots,f-1$, together with
$u_\triv$ if $\chi$ is trivial and $u_\cyc$
 if $\chi$ is cyclotomic.
Then $B$ is a basis for $U_\chi$.
\end{theorem}
{Before giving the proof, we remark that if $p=2$, then the cyclotomic character
is trivial, so the basis $B$ includes both $u_\triv$ and $u_\cyc$ and hence consists of $f+2$
elements.}

\begpf 
Define a decreasing filtration on $U_\chi$ with $\fil^0 U_\chi = U_\chi$
and $\fil^m U_\chi$ as the image of $(U_m \otimes \Fpbar(\chi^{-1}))^{\gal(M/K)}$ 
for $m \ge 1$, where $U_m = 1 + \pi^m\CO_M$. 
Setting $\gr^m U_\chi = \fil^m U_\chi/\fil^{m+1}U_\chi$, 
we see as in the proof of Theorem~\ref{thm:slopes} that
$$\dim_{\Fpbar} \gr^m U_\chi = \dim_{\Fpbar} \gr^s H^1(G_K, \Fpbar(\chi))$$
where $s = 1 + \frac{m}{e}$ if $m \ge 1$, and $s=0$ if $m=0$.
Recall that these dimensions $d_s$ are given by Theorem~\ref{thm:slopes}.
We will prove that for each $m\ge 0$, there are $d_s$ elements
of $B \cap \fil^m U_\chi$ whose images in $\gr^m U_\chi$ are linearly
independent.   It then follows that $B$ spans $U_\chi$, which suffices
since the cardinality of $B$ coincides with the dimension of $U_\chi$.

If $m > pe/(p-1)$, then $s > 1 + \frac{p}{p-1}$, so $d_s = 0$ and there
is nothing to prove.

If $m = pe/(p-1)$, then $s = 1 + \frac{p}{p-1}$, so $d_s = 0$ unless
$\chi$ is cyclotomic in which case $d_s = 1$.  Therefore it suffices
to note that $u_\cyc$ is a non-trivial element of $\fil^{pe/(p-1)} U_\chi$.

Now suppose $\frac{e}{p-1} \le m < \frac{ep}{p-1}$, so $d_s$
is the number of $i$ such that $s = 1 + \frac{m}{e} = 1 + \frac{n_{i+1}}{p^f-1}$
and $a_{i+1} \neq p$.  For each such $i$, we have $n_i' = \frac{en_{i+1}}{p^f-1} = m$
and $\tau_i' = \tau_{i+1}$, so that 
$u_i = \varepsilon_{\pi^m}(\lambda_{\tau_{i+1},\mu}) \in \fil^m U_\chi$.
We now show that the images of these $u_i$ in $\gr^m U_\chi$ are linearly independent.
We may assume that $d_s > 0$ and hence that $m$ is not divisible by $p$.
Since $m < \frac{pe}{p-1}$, we see as in the proof of Theorem~\ref{thm:slopes}
that the natural map
$$U_m/U_{m+1} \to (M^\times/ U_{m+1}) \otimes \F_p$$
is injective, so that we may identify $\gr^m U_\chi$ with
$$(U_m/U_{m+1} \otimes \Fpbar(\chi^{-1}))^{\gal(M/K)}.$$
Since the map 
$$l\otimes \Fpbar \to U_m/U_{m+1} \otimes \Fpbar$$
induced by $\varepsilon_{\pi^m}$ is an isomorphism and
the elements $\lambda_{\tau_{i+1},\mu}$ are linearly independent 
over $\Fpbar$, it follows that so are their images in 
$\gr^m U_\chi$.

Now suppose that $0 < m < \frac{e}{p-1}$.  In this case $d_s$ is the number
of $i$ such that $s = 1 + \frac{m}{e} = \frac{n_{j+1}}{p^f-1}$, $a_{j+1} \neq p-1$,
and $(a_{i+1},\ldots,a_j) = (p,p-1,\ldots,p-1)$ for some $j>i$.
For each such $i$, we have $n_i' =  \frac{en_{j+1}}{p^f-1} - e= m$
and $\tau_i' = \tau_{j+1}$, so that
$u_i = \varepsilon_{\pi^m}(\lambda_{\tau_{j+1},\mu}) \in \fil^m U_\chi$.
Note also that for distinct $i$, the corresponding $j$ are distinct mod $f$.
The proof that the images of $u_i$ in $\gr^m U_\chi$ are linearly independent
is then the same as in the preceding case.

Finally note that if $m=0$, then $s=0$, so $d_s =0$ unless $\chi$ is trivial, 
in which case $d_s = 1$.   Therefore it suffices to note that $u_\triv$ 
is not in $\fil^1 U_\chi$.
\epf

We can now define a basis for $H^1(G_K, \Fpbar(\chi))$ as the dual basis to
the one in Theorem~\ref{prop:basis}, denoting the corresponding
cohomology classes $c_\tau$ for $\tau:k \to \Fpbar$, together with $c_\nr$
if $\chi$ is trivial and $c_\tr$ if $\chi$ is cyclotomic.  We record the
construction as follows:

\begin{corollary}  \label{thm:basis}
The set consisting of the classes $c_\tau$ for $\tau \in T$, together with
$c_\nr$ if $\chi$ is trivial and $c_\tr$ if $\chi$ is cyclotomic,
forms a basis for $H^1(G_K,\Fpbar(\chi))$.
\end{corollary}

\section{Dependent pairs and admissible subsets}  \label{subsec:dependence}

We now determine the extent to which the basis for $H^1(G_K,\Fpbar(\chi))$
just constructed is independent (up to scalars) of the choices made.  We
maintain the notation of Section~\ref{subsec:basis}, so $K$ is an unramified
extension of $\Q_p$ of degree $f$ with residue field $k$, $T = \{\tau_0,\ldots,\tau_{f-1}\}$
(where $\tau_i = \tau_0 \circ\frob^i$) is the set of embeddings $k \to \Fpbar$,
and we fix a character $\chi:G_K \to \Fpbar^\times$ and write
$\chi|_{I_K} = \prod_{i=0}^f \omega_{\tau_i}^{a_i}$
where $\omega_{\tau_i}:I_K \to \Fpbar^\times$ corresponds to $\tau_i$ by
local {class} field theory and $(a_0,\ldots,a_{f-1})$ is the tame signature of $\chi$.

\subsection{Dependent pairs}
 \label{DepPairs} 
 Recall that we chose an extension $M$ of $K$ with a uniformizer $\pi$ such
that $\chi|_{G_M}$ is trivial and $M=L(\pi)$ where $L/K$ is unramified of degree
prime to $p$,  $e = [M:L]$ divides $p^f-1$, and  $\pi^e \in K^\times$. 
We introduce the following notion in order to explain how the
basis of Corollary~\ref{thm:basis} depends on the choice of $M$ and $\pi$:
\begin{definition}\label{def:dependent}
For $i,t \in \Z$ with $1 \le t \le f-1$, we say that $([i],[i+t]) \in (\Z/f\Z)^2$ is
a {\em dependent pair} if $a_{i+1} = p$, $a_{i+t+1} \neq p$, and
$$a_{i+2} = \ldots = a_{i+s} =p-1,\qquad a_{i+s+1} = \ldots = a_{i+t} = p$$
for some $s \in 1,\ldots, t$. 
\end{definition}

Note that the first (resp.~second) displayed  {chain of equations} automatically holds
if $s=1$ (resp.~$s=t$).  Note that if $a_{i+1} \neq p$ then there are no
dependent pairs of the form $([i],[j])$, and that if $a_{i+1} = p$ then the number of
dependent pairs $([i],[j])$ is either $s$ or $s-1$ where $s \in \{1,\ldots,f\}$ is such that
$$a_{i+2} = \ldots = a_{i+s} =p-1,\qquad a_{i+s+1} \neq p-1.$$
More precisely, the number of such dependent pairs is $s$ unless
$$(a_{i+2},\ldots,a_{i+s},a_{i+s+1},\ldots,a_{i+f+1}) = (p-1,\ldots,p-1,p,\ldots,p),$$
in which case it is $s-1$.  Note that there are no dependent pairs at all
if $\chi|_{I_K}$ is trivial (in which case all $a_i = p-1$) or cyclotomic
(in which case all $a_i = 1$).

\subsection{Dependence of the basis on the choice of uniformizer}
Recall that we defined a basis for 
$ U_\chi = \left( M^\times \otimes \Fpbar(\chi^{-1}) \right)^{\gal(M/K)}$
using the elements $u_i$ (for $i=0,\ldots,f-1$), $u_\triv$ (if $\chi$ is trivial)
and $u_\cyc$  (if $\chi$ is cyclotomic) defined
by equations (\ref{eqn:units}), (\ref{eqn:triv}), (\ref{eqn:cyc}).
Suppose now that we choose another uniformizer $\pi'$ (for the same $M$)
such that $(\pi')^e \in K^\times$ and accordingly define elements $u_i'$
for $i=0,\ldots,f-1$, and $u_\triv'$ if $\chi$ is trivial.
(Note that $u_\cyc$ does not depend on the choice of uniformizer.)

\begin{proposition}   \label{prop:indep_units}
For $i=0,\ldots,f-1$, the element $u_i'$ differs from a non-zero multiple of $u_i$ 
by an element of the span of
$$\{\, u_j \,|\, \mbox{$([i],[j])$ is a dependent pair}\,\}$$
and $u_\cyc$ if $\chi$ is cyclotomic.
\end{proposition}
\begpf  Recall that we require $\pi^e$ and $(\pi')^e$ to be in $K$, so
setting $\alpha = \pi'/\pi$ and $a = \ol{\alpha} \in l$, we have
$\alpha^e \in \CO_K^\times$ and $a^e \in k^\times$.

Suppose first that $\alpha = [a]$.  Note
that $\ol{\omega}_{\pi'}= \ol{\omega}_\pi \omega_a$, where $\omega_a$
is the unramified character of $\gal(L/K) \cong \gal(l/k)$
sending $g$ to $g(a)/a \in \mu_e(k)$.   Writing
$$\chi = \mu (\tau'_i \circ \overline{\omega}_\pi)^{n_i'}
 = \mu' (\tau'_i \circ \overline{\omega}_{\pi'})^{n_i'},$$
 we see $\mu = \mu' (\tau'_i \circ \omega_a)^{n_i'}$.
 Recall that $u_i$ and $u'_i$ are defined by
 $$u_i = \varepsilon_{\pi^{n_i'}}(\lambda_{\tau_i',\mu})\quad\mbox{and}\quad
     u'_i = \varepsilon_{(\pi')^{n_i'}}(\lambda_{\tau_i',\mu'}),$$
 where $\lambda_{\tau_i',\mu}$ and $\lambda_{\tau_i',\mu'}$ are
 any non-zero vectors in the corresponding eigenspaces
 $\Lambda_{\tau_i',\mu}$ and  $\Lambda_{\tau_i',\mu'}$.
 Note however that $\Lambda_{\tau_i',\mu} = (a^{n_i'}\otimes 1)\Lambda_{\tau_i',\mu'}$,
 so we may choose $\lambda_{\tau_i',\mu} = (a^{n_i'}\otimes 1)\lambda_{\tau_i',\mu'}$,
 which gives $u_i = u_i'$.   
 
 The preceding paragraph shows that we may replace $\pi$ by $[a]\pi$ and
 hence assume that $\alpha \equiv 1 \bmod \pi\CO_M$.  
 Note that $\ol{\omega}_{\pi'} = \ol{\omega}_\pi$, so $\mu = \mu'$
 and we may use the same $\lambda_{\tau_i',\mu}$ in the definitions of $u_i$ and $u_i'$.
 Since $\alpha^e \in 1 + p\CO_K$, we see that in fact $\alpha \in 1 + p\CO_K$,
 so that  $(\pi')^{n_i'} = \pi^{n_i'}\alpha^{n_i'}$ where 
 $\alpha^{n_i'} = 1 + p\delta$ for some $\delta \in \CO_K$.
 
 We now apply Lemma~\ref{lem:delta} with $x$ evaluated at $\pi^{n_i'}$.
 First note that if $a_{i+1} \neq p$, then $n_i' \ge e/(p-1)$, so $E([a]p\delta(\pi^{n_i'})^{p^m})
 \equiv 1 \bmod \pi^{ep/(p-1)}$ for all $a\in l$ and $m\ge 0$.   The lemma then implies that
 $u_i' - u_i \in \fil^{ep/(p-1)}U_\chi$, so $u_i = u_i'$ unless $\chi$ is cyclotomic,
 in which case $u_i' - u_i$ is in the span of $u_\cyc$.  
 
 Now suppose that $a_{i+1} = p$ and let $s$ be the least positive integer such that
$a_{i+s+1} \neq p-1$.    We then have $n_i' = \frac{en_{i+s+1}}{p^f-1} - e$, where
$$n_{i+s+1} = a_{i+s+1} + a_{i+s+2}p + \cdots + a_ip^{f-s-1} + p^f.$$
For $m=0,\ldots,s-1$, we see that 
$$p^m(n_{i+s+1} - (p^f-1)) = n_{i+s-m+1} - (p^f-1).$$
Since $n_{i+2} \ge 1 + p + \cdots p^{f-2} + p^f$, we see also that
$p^s(n_{i+s+1} -(p^f-1))  = p(n_{i+2} - (p^f-1)) > \frac{p^f-1}{p-1}$,
and hence that $p^sn_i' > \frac{e}{p-1}$.
Therefore $E([a]p\delta(\pi^{n_i'})^{p^m})$ is in $(\CO_M^\times)^p$ for
all $a\in l$ and $m \ge s$, so Lemma~\ref{lem:delta} implies that
$$\varepsilon_{(\pi')^{n_i'}} = 
 \varepsilon_{\pi^{n_i'}} + \sum_{m=0}^{s-1}   \varepsilon_{p\delta\pi^{n_i' p^m}}\circ \frob^m.$$
  Note also that $p\delta\pi^{n_i'p^m} = \beta\pi^{n_i'p^m + e}$
 for some $\beta \in \CO_K$, and another application of Lemma~\ref{lem:delta}
 shows that $\varepsilon_{\beta\pi^{n_i' p^m+e}} = \varepsilon_{[\ol{\beta}]\pi^{n_i' p^m+e}}$.
 It follows that 
 $$ u_i' - u_i = \tau_{i+s+1}(\ol{\beta}) \sum_{m=0}^{s-1}
 \varepsilon_{\pi^{n_i' p^m+e}}(\frob^m(\lambda_{\tau_{i+s+1},\mu})).$$
 We will show that each term is a multiple of a vector of the form $u_{i+t}$,
 where either $t \in \{1,\ldots, s-1\}$, or $t$ is the least integer such
 that $t \ge s$ and $a_{i+t+1} \neq p$.  Note that $([i],[i+t])$ is
 a dependent pair for each such $t$ (including $t = f+1$ in the case
 $(a_{i+2},\ldots,a_{i+s},a_{i+s+1},\ldots,a_{i+f+1}) = (p-1,\ldots,p-1,p,\ldots,p)$).

First consider the term with $m=0$, and note that
$n_i' + e = \frac{en_{i+s+1}}{p^f-1}$.   If $a_{i+s+1} \neq p$,
then $n_{i+s}' = \frac{en_{i+s+1}}{p^f-1} $ and  $\tau_{i+s}'  = \tau_{i+s+1}$,
so $\varepsilon_{\pi^{n_i' + e}}(\lambda_{\tau_{i+s+1},\mu}) = u_{i+s}$.
On the other hand, if $a_{i+s+1} = p$, then $n_{i+s+1}$ is divisible by $p$
and $\frac{en_{i+s+1}}{p(p^f-1)} > \frac{e}{p(p-1)}$, so Lemma~\ref{lem:shift}
implies that 
$$\varepsilon_{\pi^{en_{i+s+1}/(p^f-1)}}  = 
   \varepsilon_{-p\pi^{en_{i+s+1}/p(p^f-1)}}\circ \frob^{-1}.$$
Writing $-p = \gamma \pi^e$ for some $\gamma \in \CO_K^\times$
and noting that $n_{i+s+2} =  \frac{n_{i+s+1}}{p} + p^f-1$, we see that
this is the same as 
$\varepsilon_{\gamma \pi^{en_{i+s+2}/(p^f-1)}} \circ\frob^{-1}$,
and another application of Lemma~\ref{lem:delta} shows that we
may replace $\gamma$ by $[\ol{\gamma}]$.  Since $\frob^{-1}$ sends
$\Lambda_{\tau_{i+s+1},\mu}$ to $\Lambda_{\tau_{i+s+2},\mu}$,
we conclude that $\varepsilon_{\pi^{n_i' + e}}(\lambda_{\tau_{i+s+1},\mu})$
is a scalar multiple of $\varepsilon_{\pi^{en_{i+s+2}/(p^f-1)}}(\lambda_{\tau_{i+s+2},\mu})$.
If $a_{i+s+2} \neq p$, then this is $u_{i+s+1}$.  If $a_{i+s+2} = p$, then we
may iterate the argument to conclude that 
$\varepsilon_{\pi^{n_i' + e}}(\lambda_{\tau_{i+s+1},\mu})$ is a multiple of
$u_{i+t}$, where $t$ is the least integer such
 that $t \ge s$ and $a_{i+t+1} \neq p$.

Finally for $m=1,\ldots,s-1$, we have
 $a_{i+s-m+1} = p-1 \neq p$, so 
 $$n'_{i+s-m} = \frac{en_{i+s-m+1}}{p^f-1} = \frac{e}{p^f-1}(p^m(n_{i+s+1} - (p^f-1)) + p^f-1)
       = n'_i p^m + e$$
and  $\tau'_{i+s-m} = \tau_{i+s-m+1}$.  Since $\frob^m$ sends $\Lambda_{\tau_{i+s+1},\mu}$
 to $\Lambda_{\tau_{i+s-m+1},\mu}$, we conclude that
 $\varepsilon_{\pi^{n_i' p^m+e}}(\frob^m(\lambda_{\tau_{i+s+1},\mu}))$
 is a multiple of $u_{i+t}$ where $t = s - m$.
\epf

\subsection{Dependence of the dual basis on the choice of $M$} 
Recall 
that we defined a basis for $H^1(G_K, \Fpbar(\chi))$ as the dual basis to
the one constructed for $U_\chi$, denoting the corresponding
cohomology classes $c_\tau$ for $\tau:k \to \Fpbar$, together with $c_\nr$
if $\chi$ is trivial and $c_\tr$ if $\chi$ is cyclotomic. 
Suppose now that we make another choice of $M'$ and $\pi'$ of
the required form and denote the corresponding basis elements
$c_\tau'$ for $\tau:k \to \Fpbar$, and $c_\tr'$ if $\chi$ is cyclotomic. 
(Note that if $\chi$ is trivial, then $c_\nr$ spans $H^1_\nr(G_K,\Fpbar)$,
so, up to scalar, is independent of the choices.)

\begin{proposition} \label{prop:indep_classes}
For $i=0,\ldots,f-1$, the element $c_{\tau_i}'$ differs from a non-zero multiple
of $c_{\tau_i}$ by an element of the span of
$$\{\, c_{\tau_j} \,|\, \mbox{$([j],[i])$ is a dependent pair}\,\}$$
and $c_\nr$ if $\chi$ is trivial.
\end{proposition}
\begpf  Suppose first that the $c_\tau$ and $c_\tau'$ are defined 
using the same field $M$, but different choices of uniformizers $\pi$
and $\pi'$.   Suppose also that $\chi$ is not trivial or cyclotomic.
Define $T = (t_{ij}) \in \GL_f(\Fpbar)$ by $u_i = \sum_{i=0}^{f-1} t_{ij}u_j'$
for $i = 0,\ldots,f-1$, so that $c_{\tau_i}' = \sum_{j=0}^{f-1} t_{ji} c_{\tau_j}'$
for $i = 0,\ldots,f-1$.  The conclusion is then immediate from
Proposition~\ref{prop:indep_units}, which shows that
$t_{ii} \neq 0$ for each $i$, and that
$t_{ij} = 0$ unless $i=j$ or $([i],[j])$ is a dependent pair.
If $\chi$ is trivial or cyclotomic, then there are no dependent
pairs, and the conclusion is again immediate from
Proposition~\ref{prop:indep_units}.

Now suppose that $M$ and $M'$ are any two extensions of
$K$ of the required form.  By symmetry, we may replace $M'$
with a larger extension satisfying the hypotheses, and hence
assume that if $M = L(\pi)$, then $M' = L'(\pi')$ where $L'$
is an unramified extension of $L$ of degree prime to $p$
and $(\pi')^d = \pi$ where $de$ divides $p^f-1$.
By the preceding paragraph, we may assume that the
$c_\tau'$ are defined using the uniformizer $\pi'$.

Note that we have used the isomorphisms of class field theory in order to identify
$H^1(G_K,\Fpbar(\chi))$ with both $\hom_{\Fpbar}(U_\chi,\Fpbar)$ and
with $\hom_{\Fpbar}(U_\chi',\Fpbar)$, where
$$ U_\chi = \left( M^\times \otimes \Fpbar(\chi^{-1}) \right)^{\gal(M/K)}
\quad\mbox{and}\quad
U_\chi' = \left( (M')^\times \otimes \Fpbar(\chi^{-1}) \right)^{\gal(M'/K)}.$$
Recall that this identification is compatible with the isomorphism
$U_\chi' \to U_\chi$ induced by the norm map from $(M')^\times$
 to $M^\times$.  Denoting this isomorphism $\nu_{M'/M}$ and the basis
 elements for $U_\chi'$ by $u_i'$, it
 suffices to prove that $\nu_{M'/M}(u'_i)$ is a multiple of $u_i$ for
 $i=0,\ldots,f-1$, and similarly for $u_\cyc'$ and $u_\cyc$ if $\chi$ is cyclotomic.
 
 With our choices of $\pi$ and $\pi'$, the map $\varepsilon_{\pi^{n_i'}}$
 appearing in the definition of $u_i$ is simply the restriction to $l \otimes \Fpbar$
 of the one in the definition of $u_i'$.  Note also that the embeddings $\tau_i'$
 and unramified characters $\mu$ are the same for $M$ and $M'$.
 Therefore
 $$\nu_{M'/M}(u_i') = \sum_{g\in\gal(M'/M)} \varepsilon_{\pi^{n_i'}}\lambda'_{\tau_i',\mu}
   = d \varepsilon_{\pi^{n_i'}} \tr_{l'/l}(\lambda'_{\tau_i',\mu}),$$
where $\lambda'_{\tau_i',\mu}$ is in the $\mu$-eigenspace for the action of
$\gal(l'/k)$ on $l'\otimes_{k,\tau_i'}\Fpbar$.    
The conclusion follows from the observation that 
$\tr_{l'/l}(\lambda'_{\tau_i',\mu}) \in \Lambda_{\tau_i',\mu}$.
Finally, if $\chi$ is cyclotomic, then the argument for $u_\cyc'$ is similar.
\epf

\subsection{Admissible sets}

\begin{definition}\label{def:admissible}
We say that a subset $J \subset \Z/f\Z$ is {\em admissible} if for all
dependent pairs $([j],[i])$, we have that if $[i] \in J$, then $[j] \in J$.
We say that a subset $J \subset T$ is {\em admissible} if the corresponding
subset of $\Z/f\Z$ is admissible.
\end{definition}

The following is immediate from Proposition~\ref{prop:indep_classes}:

\begin{corollary}  \label{cor:indep_spaces} If $J \subset T$ is admissible,
then the span of the set $\{\,c_\tau\,|\,\tau \in J\,\}$ in
$H^1(G_K,\Fpbar(\chi))/H^1_\ur(G_K,\Fpbar(\chi))$
is well-defined, i.e., independent of the choice of
$M$ and $\pi$.
\end{corollary}

Finally we give some criteria for admissibility in terms of the subspaces
of $H^1(G_K,\Fpbar(\chi))$ which were defined in Section~\ref{subsec:fil}
using the ramification filtration.  Note that since $c_\nr \in H^1_\nr(G_K,\Fpbar(\chi))$
(if $\chi$ is trivial) and $c_\tr \not\in H^1_\ar(G_K,\Fpbar(\chi))$
(if $\chi$ is cyclotomic), we always have that 
$\{\,c_\tau\,|\,\tau \in T\,\}$ is a basis for the $f$-dimensional space
$H^1_\ar(G_K,\Fpbar(\chi))/H^1_\nr(G_K,\Fpbar(\chi))$.

\begin{theorem} \label{thm:admissibility}
With the above notation we have:  
\begin{enumerate} 
\item If $\tau \in T$, then the following hold:
\begin{enumerate}
\item $\{\tau\}$ is admissible if and only if $c_\tau \in   H^1_\cg(G_K,\Fpbar(\chi))$;
\item $T - \{\tau\}$ is admissible if and only if $c_\tau \not\in  H^1_\mr(G_K,\Fpbar(\chi))$.
\end{enumerate}
\item The following are equivalent:
\begin{enumerate}
\item $\chi$ is generic;
\item all subsets of $T$ are admissible;
\item $H^1_\cg(G_K,\Fpbar(\chi)) = H^1_\ar(G_K,\Fpbar(\chi))$;
\item $H^1_\mr(G_K,\Fpbar(\chi)) = H^1_\nr(G_K,\Fpbar(\chi))$.
\end{enumerate}
\end{enumerate}
\end{theorem}
\begpf To prove part 1), let $\tau = \tau_i$.  
From the proof of Theorem~\ref{prop:basis},
we see that $c_\tau$  {is in} $\fil^s H^1(G_K,\Fpbar(\chi))$ but not in
$\fil^{<s} H^1(G_K,\Fpbar(\chi))$, where $s = 1 + \frac{n_i'}{e}$
and $n_i'$ is as in the definition of the classes $u_i$.

For 1a), note that  $\{ i \}$ fails to be admissible if and only if 
$a_{i+1} \neq p$ and $(a_j,\ldots,a_i) = (p,p-1,\ldots,p-1)$ for some $j$
with $i-f+1 < j \le i$, which in turn is equivalent to $n_i' > e$.
Therefore $\{\tau\}$ is admissible if and only if $c_\tau \in 
\fil^2 H^1(G_K,\Fpbar(\chi)) = H^1_\cg(G_K,\Fpbar(\chi))$.

For 1b), note that $T - \{\tau\}$ is admissible if and only if
$a_{i+1} \neq p$ or $\chi|_{I_K}$ is not cyclotomic, which in turn is
equivalent to $n_i' \ge e/(p-1)$.   Therefore $T - \{\tau\}$ is admissible
if and only if 
$ {c_\tau \not \in \fil^{<\frac{p}{p-1}}} H^1(G_K,\Fpbar(\chi))
 = H^1_\mr(G_K,\Fpbar(\chi))$.
 
Turning to part 2), note that 2a) and 2b) are both equivalent to the
condition that there be no dependent pairs, which in turn is equivalent 
to the admissibility of all singletons.  The equivalence of 2b) and 2c) thus follows
from  {1a)} and the fact that the $\{\,c_\tau\,|\,\tau \in T\,\}$ span
$H^1_\ar(G_K,\Fpbar(\chi))/H^1_\nr(G_K,\Fpbar(\chi))$. The
equivalence of 2a) and 2d) is immediate from part 5) of
Corollary~\ref{cor:slopes}.
\epf

\section{Distinguished subspaces}  \label{subsec:subspaces}
We now return to the setting of Section~\ref{subsec:weights}, so $K$ is an
unramified extension of $\Q_p$ of degree $f$ with residue field $k$, $T$
is the set of embeddings $k \to \Fpbar$, and $\rho:G_K \to \GL_2(\Fpbar)$
is a continuous representation.  We assume further that $\rho$ is reducible,
so $ {\rho} \sim \begin{pmatrix}\chi_1&*\\ 0&\chi_2\end{pmatrix}$
for some characters $\chi_1$ and $\chi_2$ of $G_K$.
We let $\chi = \chi_1\chi_2^{-1}$ and let $c_\rho \in H^1(G_K,\Fpbar(\chi))$ denote the
extension class associated to $\rho$.  

Recall that we have defined a set $W'(\chi_1,\chi_2)$ of certain pairs
$(V,J)$, where $V$ is a Serre weight (i.e., an irreducible representation
of $\GL_2(k)$ over $\Fpbar$) and $J$ is a subset of $T$ and for each
$(V,J) \in W'(\chi_1,\chi_2)$ a certain subspace $L_{V,J}$ of 
$H^1(G_K,\Fpbar(\chi))$.   We then define {d} $L_V$ as the union of
the subspaces $L_{V,J}$ such that $(V,J) \in W'(\chi_1,\chi_2)$,
and $W(\rho)$ as the set of $V$ such that $c_\rho \in {L_V}$ {fixed typo}.
In this section we give a conjectural description of $L_V$ in
terms of the basis for $H^1(G_K,\Fpbar(\chi))$ of Corollary~\ref{thm:basis}.

As in the preceding sections, we choose an embedding $\tau_0 \in T$,
set $\tau_i = \tau_0\circ \frob^i$ (for $i\in \Z$ or $\Z/f\Z$), and let
$\omega_{\tau_i} : I_K \to \Fpbar^\times$ denote the corresponding
fundamental character.   We let $\chi = \chi_1\chi_2^{-1}$
and write
$$\chi|_{I_K} = \prod_{i = 0}^{f-1}  \omega_{\tau_i}^{a_i},$$
where $(a_0,\ldots,a_{f-1})$ is the tame signature of $\chi$. 
We will often interchange $\tau_i$ and $i$
in the notation, and thus identify  $T$ with $\Z/f\Z$ (except in the notation
for fundamental characters where this could lead to confusion).

\subsection{Shifting functions $\delta$ and $\mu$} For any subset $J$ of $\Z/f\Z$, we will define a subset $\mu(J)$ of $\Z/f\Z$.
First we define a function
$\delta: \Z \to \Z$ depending on the integers $a_i$ as follows: 
If $j \in \Z$ we let $\delta(j) = j$ unless
$$(a_{i+1},a_{i+2},\ldots,a_j) = (p,p-1,\ldots,p-1)$$
for some $i < j$ (necessarily unique), in which case we let $\delta(j) = i$.
Note that $\delta$ induces a function $\Z/f\Z \to \Z/f\Z$, which we also
denote by $\delta$.
If $\delta(J) \subset J$, then we let $\mu(J) = J$. Otherwise we choose
some $[i_1] \in \delta(J) \setminus J$ and let $j_1$ be the least integer
such that $j_1 > i_1$, $[j_1] \in J$ and $\delta(j_1) = i_1$.  Now write
$J = \{[j_1],\ldots,[j_r]\}$ with $j_1 < j_2 < \cdots < j_r < j_1 + f$,
and define $i_{\kappa}$  for $\kappa=2,\ldots,r$ inductively as follows:
$$i_\kappa = \left\{\begin{array}{cl} \delta(j_\kappa),&\mbox{if $i_{\kappa-1} < \delta(j_\kappa)$,}\\
j_\kappa ,&\mbox{otherwise.}\end{array}\right.$$
We then have $i_1 < i_2 < \cdots < i_r < i_1 + f$, and we set
$\mu(J) = \{[i_1],\ldots,[i_r]\}$.  One checks easily that this is
independent of the choice of $i_1$.  Note that by construction
we have $\delta(J) \subset \mu(J) \subset \delta(J) \cup J$.

\begin{lemma} \label{lem:admissible}
The set $\mu(J)$ is admissible.
\end{lemma}
\begpf  Suppose that ${([i+t],[i])}$ is a dependent pair with
notation as in Definition~\ref{def:dependent}.  We must show that
if $[i+t] \in \mu(J)$, then $[i] \in \mu(J)$.  {Recall from \S\ref{DepPairs} that $s$ is such that $a_{i+s}=p-1$ but $a_{i+s+1}=p$.}  Note that 
$$\delta(i+1) = \cdots = \delta(i+s) =i, \quad\mbox{and}\quad
\delta(i+\nu) = i+\nu - 1 \mbox{\ for $\nu = s+1,\ldots,t$.}$$
In particular $\delta([i+t]) \neq [i+t]$, and since $a_{i+t+1} \neq p$,
it follows that $[i+t]$ is not in the image of $\delta$.  
If $[i+t] \in \mu(J)$, we must therefore have $i+t = j_\kappa = i_\kappa$
for some choice of $i_1$ and some $\kappa \in \{2,\ldots,r\}$.
If $s < t$, then the resulting inequalities
$$j_\kappa - 1 = \delta(j_\kappa) \le i_{\kappa - 1} \le j_{\kappa - 1}
        < j_\kappa$$
imply that $i_{\kappa - 1} = j_{\kappa - 1} = i+ t - 1$ and $\kappa \ge 3$. 
Repeating the argument shows that for $\nu =2,\ldots,s-t$, we have
$i_{\kappa - \nu} = j_{\kappa - \nu} = i + t - \nu$ and $\kappa \ge \nu + 2$.
In particular $[i+s] \in J$, and hence $[i] = \delta([i+s]) \in \mu(J)$.
\epf

\subsection{Explicit distinguished subspaces} Now let $V = V_{\vec{d},\vec{b}}$ and suppose that $(V,J) \in W'(\chi_1,\chi_2)$
for some $J \subset T$.   Then there is a unique $J_\max \subset \Z/f\Z$ such that 
$(V,J_\max) \in W'(\chi_1,\chi_2)$, and $J_\max$ satisfies
the two conditions:
\begin{itemize}
\item if $(b_i,b_{i+1},\ldots, b_{j-1},b_j) = (p,p-1,\ldots,p-1,1)$ for some $i < j$
such that  $i,i+1,\ldots,j-1 \not\in J_\max$, then $j\not\in J_\max$;
\item if $(b_0,\ldots,b_{f-1}) = (p-1,p-1,\ldots,p-1)$, or if $p=2$ and 
$(b_0,\ldots,b_{f-1}) = (2,2,\ldots,2)$, then $J_\max \neq \emptyset$.
\end{itemize}
(This is proved in \cite{gls} for $p>2$, but one easily checks that it holds also
for $p=2$.)
We then define
$L^\AH_V   {\subset} H^1(G_K,\Fpbar(\chi))$ to be the span of
$\{\,c_\tau\,|\,\tau \in \mu(J_\max)\,\}$ together with $c_\nr$ if
$\chi$ is trivial, unless $\chi$ is cyclotomic
and $V = V_{\vec{d},\vec{b}}$ with 
$J_\max = T$
and $(b_0,\ldots,b_{f-1}) = (p,p,\ldots,p)$, in which
case $L^\AH_V = H^1(G_K,\Fpbar(\chi))$ (i.e., we include
$c_\tr$ as well).   By Corollary~\ref{cor:indep_spaces}
and Lemma~\ref{lem:admissible} the space $L^\AH_V$
is well-defined, i.e., independent of the choices made in
Section~\ref{subsec:basis}.  (The superscript $\AH$ refers to
the use of the Artin--Hasse exponential in its definition.)

We now state our conjectural explicit description of the subspaces appearing in the
recipe for the weight:
\begin{conjecture}  \label{conj:spaces}
If $(V,J) \in W'(\chi_1,\chi_2)$ for some $J$, then $L_V = L^\AH_V$.
 \end{conjecture}
 
Recall that \cite{gls} proves that if $p>2$ and $(V,J) \in W'(\chi_1,\chi_2)$,
then $L_{V,J} \subset L_{V,J_\max}$, so that $L_V$ can be replaced
by $L_{V,J_\max}$ in the statement of Conjecture~\ref{conj:spaces}
if $p>2$.  Since $L_V^\AH$ is a subspace of $H^1(G_K,\Fpbar(\chi))$
of dimension at most that of $L_{V,J_\max}$, Conjecture~\ref{conj:spaces}
implies that $L^\AH_V = L_{V,J_\max}$, and hence the assertion that
$L_{V,J} \subset L_{V,J_\max}$ for all $(V,J) \in W'(\chi_1,\chi_2)$
still holds for $p=2$.

\subsection{A weight-explicit Serre's conjecture}
Finally we record a more explicit form of Conjecture~\ref{conj:serre1}.
For $\rho$ as above, define $W^\AH(\rho)$ to be the set of $V$ such
that $(V,J) \in W'(\chi_1,\chi_2)$ for some $J$ and $c_\rho \in L_V^\AH$.
For irreducible $\rho:G_K \to \GL_2(\Fpbar)$, define $W^\AH(\rho) = W(\rho)$.

Suppose now that $\rho:G_F \to \GL_2(\Fpbar)$ is continuous, irreducible
and totally odd.    Combining Conjectures~\ref{conj:serre1} and~\ref{conj:spaces}
then yields:

\begin{conjecture} \label{conj:serre2}  {The representation}
$\rho$ is modular of weight $V  = \otimes_{\{\gp\in S_p\}} V_\gp$
if and only if $V_\gp \in W^\AH(\rho|_{G_{F_\gp}})$ for all $\gp \in S_p$.
\end{conjecture}

If $\rho$ is modular (of some weight) and satisfies the hypotheses under
which the weight part of Serre's Conjecture is known (by \cite{gk} or \cite{newton}),
then Conjecture~\ref{conj:serre2} is immediate from Conjecture~\ref{conj:spaces}.

\section{The quadratic case} \label{sec:quadratic}
In this section we delineate the possibilities for the spaces of extensions
$L_V^\AH$ and the sets of Serre weights $W^\AH(\rho)$ in the case $f=2$.
We refer the reader to the forthcoming paper \cite{ddr} for a more detailed
discussion of the situation for arbitrary $f$ and the underlying combinatorics.

Suppose now that $K$ is the unramified quadratic extension of $\Q_p$,
and $\rho: G_K \to \GL_2(\Fpbar)$ is a continuous representation.

\subsection{Three reducible cases}
Suppose first that $\rho$ is reducible, so that 
$\rho \sim \begin{pmatrix}\chi_1&*\\0&\chi_2\end{pmatrix}$
for some characters $\chi_1,\chi_2:G_K \to \Fpbar^\times$.
Twisting
by $\chi_2^{-1}$, we may assume $\chi_2 = 1$, and we write
$\chi$ for $\chi_1$ and $c_\rho \in H^1(G_K,\Fpbar(\chi))$
for the associated extension class.
Choosing an embedding $\tau_0:k\to \Fpbar$, we may write 
$\chi|_{I_K}=\omega_{\tau_0}^a$ with $p+1 \le a < p^2 + p$.
We let $(a_0,a_1)$ denote the tame signature of $\chi$, so
$a= a_0 + a_1p$;   altering our choice of $\tau_0$, we may further assume that
$1\le a_0  \le p -1$ and $a_0 \le a_1 \le p$.
We now divide our analysis into three cases,  %
following 
 the terminology introduced after Theorem~\ref{thm:slopes}:
\begin{itemize}
\item[I)] $\chi$ is primitive and generic: $1 \le a_0 < a_1 < p$;
\item[II)]  $\chi$ is imprimitive and generic: $1 \le a_0 = a_1 < p$;
\item[III)]  $\chi$ is primitive and non-generic: $1 \le a_0 < a_1 = p$.
\end{itemize}
Thus the analysis in Case I is simplest, and the other two cases represent the
two main complications that can occur.   Note that Case II occurs precisely when
$\chi$ has absolute niveau $1$, so $\chi|_{I_K} = \omega^c$ where $\omega$ is
the cyclotomic character and $1 \le c \le p-1$,
and Case III occurs precisely when $\chi|_{I_K} = \omega_\tau^c$
for some $\tau$ and $c$ with $1 \le c \le p-1$.  Note also for $f=2$
(or indeed any prime $f$), $\chi$ cannot be both imprimitive and non-generic.

\subsection{Case I} In Case I, the elements of $W'(\chi,1)$ are the pairs $(V_{\vec{d},\vec{b}},J)$
given by the columns of the table:
\begin{equation}\label{eqn:Wtable1}
\begin{tabular}{  >{$}c<{$}|  >{$}c<{$}  >{$}c<{$}  >{$}c<{$}  >{$}c<{$}  }\toprule
J                   &                T            &  \{0\}               &    \{1\}                 & \emptyset  \\ \midrule
\vec{d}         &           (0,0)          & (p-1,a_1-1)    & (a_0-1,p-1)      & (a_0,a_1) \\
\vec{b}         &   (a_0,a_1)         & (a_0+1,p-a_1)&(p-a_0,a_1+1)&(b_0,b_1) \\
\bottomrule
\end{tabular}
\end{equation}
where 
$$(b_0,b_1) = \left\{ \begin{array}{ll}
(p-1-a_0,p-1-a_1), & \mbox{if $a_1 < p-1$,}\\
(p-2-a_0,p),& \mbox{if $a_1 = p-1$ and $a_0 < p-2$,}\\
(p,p-1),&\mbox{if $(a_0,a_1) =  (p-2,p-1)$.}
\end{array}\right.$$
For each $V = V_{\vec{d},\vec{b}}$ in the table, there is a unique
$J \in S_V(\chi,1)$, so that $J = J_\max$; moreover $J = \mu(J)$
is admissible.  If $J=T$, then $L_V^\AH$ is the whole two-dimensional
space $H^1(G_K,\Fpbar(\chi))$, and if $J = \emptyset$, then $L_V^\AH = 0$,
but if $J= \{i\}$ for $i = 0$ or $1$, then $L_V^\AH$ is a one-dimensional
subspace of $H^1(G_K,\Fpbar(\chi))$ which we will simply denote $L_i$.
We then have four possibilities for $W^\AH(\rho)$:
\begin{itemize}
\setlength{\itemindent}{0.3in}
\item[Ia)] $\{\,V_{(0,0),(a_0,a_1)}\,\}$ if $c_\rho \not\in L_0 \cup L_1$;
\item[Ib${}_1$)] $\{\,V_{(0,0),(a_0,a_1)}, V_{(a_0-1,p-1),(p-a_0,a_1+1)}\,\}$ if $c_\rho \in  L_1$,
         $c_\rho\neq 0$;
\item[Ib${}_2$)] $\{\,V_{(0,0),(a_0,a_1)}, V_{(p-1,a_1-1),(a_0+1,p-a_1)}\,\}$ if $c_\rho \in  L_0$,
         $c_\rho\neq 0$;
\item[Ic)] $\{\,V_{(0,0),(a_0,a_1)}, V_{(p-1,a_1-1),(a_0+1,p-a_1)}, V_{(a_0-1,p-1),(p-a_0,a_1+1)},
          V_{(a_0,a_1),(b_0,b_1)}\,\}$ if $c_\rho = 0$.
\end{itemize}

We now proceed to describe the subspaces $L_0$ and $L_1$.

With notation as in Section~\ref{subsec:basis}, we have $n_0 = a = a_0 + pa_1$
and $n_1 = a_1 + pa_0$, so that $n_1 < n_0$.  Choosing a tamely ramified extension
$M$ with uniformizer $\pi$, residue field $l$, and  ramification degree $e$ as in
that section, we have
$n_0' = en_1/(p^2-1)$, $n_1' = en_0/(p^2-1)$, and
$$\chi = \mu (\tau_1\circ\ol{\omega}_\pi)^{n_0'}
     =  \mu (\tau_0\circ\ol{\omega}_\pi)^{n_1'},$$
where $\omega_\pi(g) = g(\pi)/\pi$.

Theorem~\ref{prop:basis} provides a basis $\{u_0, u_1\}$ for 
$U_\chi = (M^\times \otimes \Fpbar(\chi^{-1}))^{\gal(M/K)}$ with
$$u_i \in \varepsilon_{\pi^{n_i'}}(l \otimes \Fpbar) =
E([l ]\pi^{n_i'}) \otimes \Fpbar.$$
We can therefore describe the elements of the dual basis $\{c_0,c_1\}$ for
$$H^1(G_K,\Fpbar(\chi)) 
 \cong \hom_{\gal(M/K)}(M^\times,\Fpbar(\chi))$$
by specifying their values on the elements of $M^\times$ of the form
$E([a]\pi^{n_i'})$ for $a \in l$ and $i=0,1$.
We find that $c_0$ and $c_1$ are defined (up to scalars) by the
homomorphisms
$$\begin{array}{rcl}c_0(E([a]\pi^{n_1'})) = 0,&\qquad& \displaystyle c_0(E([a]\pi^{n_0'}) ) = 
   \sum_{g\in\gal(l/k)}  \mu^{-1}(g) \tilde{\tau}_0(ga^p)
   \\&&\\
\mbox{and}\quad c_1(E([a]\pi^{n_0'})) = 0,&\qquad  &c_1(E([a]\pi^{n_1'}))  = 
 \displaystyle   \sum_{g\in\gal(l/k)}  \mu^{-1}(g) \tilde{\tau}_1(ga^p)\end{array} $$
for any choices of embeddings $\tilde{\tau}_i: l \to \Fpbar$
extending the $\tau_i$.  Indeed it is straightforward to check
that $c_i(hx) = \chi(h) c_0(x)$ for $x = E([a]\pi^{n_i'})$, we clearly
have $c_0(u_1) = c_1(u_0) = 0$, and the following lemma
shows that $c_0$ and $c_1$ are not identically $0$.
For the lemma, we momentarily drop the assumptions that
$[k:\F_p] = 2$ and that $[l:k]$ is not divisible by $p$.

\begin{lemma} \label{lem:sillytrick}
Suppose that $k \subset l$ are finite extensions of $\F_p$,
$\mu:\gal(l/k) \to \Fpbar^\times$ is a character
and $\tilde{\tau}:l \to \Fpbar$ is an embedding.  Then the function
$f:l\to \Fpbar$ defined by
$f(a) = \sum_{g\in\gal(l/k)}  \mu(g) \tilde{\tau}(ga)$
is not identically zero.
\end{lemma}
\begpf  Suppose that $f(a) = 0$ for all $a\in l$.  Let $\F$ denote
the subfield of $\Fpbar$ generated by the values of $\mu$, and
let $r = [\F:\F_p]$.  For $i =0,\ldots,r-1$, consider the function
$f^{(i)}:l \to \Fpbar$ defined by 
$f^{(i)}(a) = \sum_{g\in\gal(l/k)}  \mu^{p^i}(g) \tilde{\tau}(ga)$.
Since $f^{(i)}(a^{p^i}) = (f(a))^{p^i} = 0$ for all
$a \in l$, the function $f^{(i)}$ is identically zero, and
therefore so is the function $h:l \to \Fpbar$ defined by
$$h(a) = \sum_{i=0}^{r-1}f^{(i)}(a)  = \sum_{g\in\gal(l/k)}  \tr_{\F/\F_p}(\mu(g)) \tilde{\tau}(ga).$$
Taking $a$ so that $\{\,ga\,|\,g\in \gal(l/k)\,\}$ is a normal basis
for $l/k$, the $\tilde{\tau}(ga)$ are linearly independent over
$\tilde{\tau}(k)$, and hence over $\F_p$.  It follows that
$\tr_{\F/\F_p}(\mu(g)) = 0$ for all $g \in \gal(l/k)$.  
Since the values $\mu(g)$ span $\F$ as a vector space
over $\F_p$, this implies that $\tr_{\F/\F_p}$ is identically zero,
yielding a contradiction.
\epf 

We thus obtain the criterion that $c_\rho \in L_i$ if and only
if $E([a]\pi^{en_i/(p^2-1)}) \in \ker(c_\rho)$ for all $a \in l$.
Since $n_1 < n_0$, this provides a description of $L_0$
in terms of the ramification filtration on cohomology
defined in Section~\ref{subsec:fil}.   By Theorem~\ref{thm:slopes},
we have 
$$\dim_{\Fpbar} \fil^s(H^1(G_K,\Fpbar(\chi))) = \left\{\begin{array}{ll}
0,&\mbox{if $s < 1 + n_1/(p^2-1)$,}\\
1,&\mbox{if $1 + n_1/(p^2-1) \le s < 1 + n_0/(p^2-1)$,}\\
2,&\mbox{if $1 + n_0/(p^2-1) \le s$.}\end{array}\right.$$
We thus see that 
\begin{equation}\label{criterion:slope1}
L_0 =  \fil^{1+n_1/(p^2-1)}(H^1(G_K,\Fpbar(\chi)))
 = \fil^{< 1+n_0/(p^2-1)}(H^1(G_K,\Fpbar(\chi))),\end{equation}
so that $c_\rho \in L_0$ if and only if $G_K^{n_0/(p^2-1)} \subset \ker(\rho)$.
The space $L_1$ cannot be described in terms of the ramification filtration,
but it can still be characterized in terms of splitting fields.  Indeed if we let
$N$ denote the splitting field over $M$ of $\rho$, then we have
\begin{equation}\label{criterion:norm} c_\rho \in L_i \quad\mbox{if and only if}\quad
E([a]\pi^{en_i/(p^2-1)}) \in \Nm_{N/M}(N^\times)
\mbox{\ for all $a\in l$}.\end{equation}

\subsection{Case II} We now turn to Case II, where the tame signature $(a_0,a_0)$ has period $1$. 
Then $W'(\chi,1)$ is given exactly as in
(\ref{eqn:Wtable1}) with the following changes:
\begin{itemize}
\item if $a_0=1$, then we also have $\vec{d} = (0,0)$, $\vec{b} = (p,p)$
    for $J = T$;
\item if $a_0= p-2$, then we also have $\vec{d} = (p-2,p-2)$, $\vec{b} = (p,p)$
    for $J = \emptyset$;
\item if $a_0=p-1$, then take $\vec{b} = (p-1,p-1)$ for $J = \emptyset$, and
we have the following additional elements:
\begin{equation*}
\begin{tabular}{  >{$}c<{$}|   >{$}c<{$}    >{$}c<{$}   >{$}c<{$}   }\toprule
J                   &                \{0\}               &    \{1\}                 & \emptyset \,\, (\text{if}\,\, p = 2) \\\midrule
\vec{d}         &           (p-2,p-1)          & (p-1,p-2)    & (0,0)\\
\vec{b}         &            (1,p)                  & (p,1)          &(2,2)\\
\bottomrule
\end{tabular}
\end{equation*}
\end{itemize}
For each $V$ we still have a unique $J \in S_V(\chi,1)$ unless
$a_0 = p-1$, in which case each $S_V(\chi,1)$ has two elements,
and the ones appearing in the last bullet above are precisely those for
which $J \neq J_\max$.   Note that every $J$ arises as $J_\max$ for
some $V$ unless $a_0 = p-1$, in which case $J_\max = \emptyset$
does not arise.  Moreover $J_\max$ uniquely determines $V$
unless $a_0 = 1$, in which case $V_{(0,0),(1,1)}$ and $V_{(0,0),(p,p)}$
both have $J_\max = T$, or $a_0 = p-2$, in which case $V_{(p-2,p-2),(0,0)}$
and $V_{(p-2,p-2),(p,p)}$ both have $J_\max = \emptyset$.
It is still the case that $J = \mu(J)$ is admissible
for every $J$.    

If $J_\max = \emptyset$, then $L_V^\AH = 0$.
If $J_\max = T$, then $L^\AH_V = H^1(G_K,\Fpbar(\chi))$ unless $\chi$
is cyclotomic and $V = V_{(0,0),(1,1)}$, in which case 
$L^\AH_V = H^1_{\pr}(G_K,\Fpbar(\chi)) = H^1_{\ar}(G_K,\Fpbar(\chi)$ 
has codimension one in $H^1(G_K,\Fpbar(\chi))$.
If $J_\max = \{i\}$ for $i = 0$ or $1$, then writing simply $L_i$ for $L_V^\AH$,
we have the sequence of inclusions of subspaces with codimension one:
\begin{equation}\label{eqn:inclusion}
H^1_\nr(G_K,\Fpbar(\chi)) \subset L_i \subset H^1_\ar(G_K,\Fpbar(\chi)).
\end{equation}

We now list the possibilities for $W^\AH(\rho)$.

If $a_0 = 1$, then $W^\AH(\rho)$ is:
\begin{itemize}
\setlength{\itemindent}{0.2in}
\item[$\mathrm{IIz)}$] $\{\, V_{(0,0),(p,p)}\,\}$ if $c_\rho \not\in H^1_\ar(G_K,\Fpbar(\chi))$;
\item[$\mathrm{IIa)}$]  $\{\, V_{(0,0),(p,p)}, V_{(0,0),(1,1)}\,\}$ if $c_\rho \in H^1_\ar(G_K,\Fpbar(\chi)) - (L_0 \cup L_1)$;
\item[$\mathrm{IIb}_1)$] $\{\, V_{(0,0),(p,p)}, V_{(0,0),(1,1)},V_{(0,p-1),(p-1,2)}\,\}$ if $c_\rho \in L_1 - L_0$;
\item[$\mathrm{IIb}_2)$] $\{\, V_{(0,0),(p,p)}, V_{(0,0),(1,1)},V_{(p-1,0),(2,p-1)}\,\}$ if $c_\rho \in L_0 - L_1$;
\item[$\mathrm{IIc})$] $\{\, V_{(0,0),(p,p)}, V_{(0,0),(1,1)},V_{(p-1,0),(2,p-1)},V_{(0,p-1),(p-1,2)},V_{(1,1),(p-2,p-2)}\,\}$\\
if $c_\rho \in L_0 \cap L_1 = H^1_\nr(G_K,\Fpbar(\chi)),$
\end{itemize}
where in Case IIc) we omit $V_{(1,1),(p-2,p-2)}$ if $p=2$ and add $V_{(1,1),(3,3)}$ if $p=3$.
Note that Case IIz) is only possible if $\chi$ is cyclotomic, and recall that 
$H^1_\nr(G_K,\Fpbar(\chi)) = 0$ unless $\chi$ is trivial  (which implies here that $p=2$).

If $2 \le a_0 \le p-1$, then $W^\AH(\rho)$ is:
\begin{enumerate}
\setlength{\itemindent}{0.2in}
\item[$\mathrm{IIa}')$] $\{\,V_{(0,0),(a_0,a_0)}\,\}$ if $c_\rho \not\in L_0 \cup L_1$;
\item[$\mathrm{IIb}_1')$] $\{\,V_{(0,0),(a_0,a_0)}, V_{(a_0-1,p-1),(p-a_0,a_0+1)}\,\}$ if $c_\rho \in  L_1 - L_0$;
\item[$\mathrm{IIb}_2')$] $\{\,V_{(0,0),(a_0,a_0)}, V_{(p-1,a_0-1),(a_0+1,p-a_0)}\,\}$ if $c_\rho \in  L_0 - L_1$;
\item[$\mathrm{IIc}')$] $\left\{\,V_{(0,0),(a_0,a_0)}, V_{(p-1,a_0-1),(a_0+1,p-a_0)}, V_{(a_0-1,p-1),(p-a_0,a_0+1)},\right.\\
\left. \qquad\qquad\qquad V_{(a_0,a_0),(p-1-a_0,p-1-a_0)}\,\right\}$ if $c_\rho \in L_0 \cap L_1 = H^1_\nr(G_K,\Fpbar(\chi))$,
\end{enumerate}
where in Case IIc$'$) we omit $V_{(1,1),(p-1-a_0,p-1-a_0)}$ if $a_0 = p-1$ and add 
$V_{(1,1),(p,p)}$ if $a_0 = p-2$.  (Recall again 
that $H^1_\nr(G_K,\Fpbar(\chi)) = 0$ unless $\chi$ is trivial, in which case $a_0 = p-1$.)

We now turn to the description of the subspaces $L_i$.  The main
difference from Case I is that we now have $n_0 = n_1 = a = a_0(1+p)$, so that
$n_0' = n_1'$ in the notation of Section~\ref{subsec:basis}.
(Note also that we may choose $e$ to divide $p-1$.)
Another difference is that $\chi$ may be trivial or cyclotomic,
so that $H^1(G_K,\Fpbar(\chi))$ and $U_\chi$ may have dimension
greater than two.  However from the inclusions~(\ref{eqn:inclusion}) we
see that it suffices to describe the image $L_i'$ of $L_i$ in the quotient
$$\begin{array}{rcl} H_\ar^1(G_K,\Fpbar(\chi))/H^1_\nr(G_K,\Fpbar(\chi))
   & \cong & \hom_{\gal(M/K)} (\CO_M^\times/U_m, \Fpbar(\chi))\\
    &\subset &\hom_{\gal(M/K)} (\CO_M^\times, \Fpbar(\chi)),\end{array}$$
where $U_m = 1 + \pi^m\CO_M$ for $m = \lceil ep/(p-1) \rceil$.
This quotient has a basis $\{c_0',c_1'\}$ where $c_i'$ spans
$L_i'$ and is determined by its values on 
elements of the form $E([a]\pi^{n_0'})$ for $a\in l$
by the formula
\begin{equation}\label{criterion:alignment} \displaystyle c_i'([a]\pi^{n_0'})  = 
   \sum_{g\in\gal(l/k)}  \mu^{-1}(g) \tilde{\tau}_i(ga^p),\end{equation}
where $\tilde{\tau}_i:l \to \Fpbar$ is any choice of embedding extending $\tau_i$.
Indeed it follows from the definitions of the elements $u_j$
that $c_i' \in L_i'$ and from Lemma~\ref{lem:sillytrick} that $c_i' \neq 0$.

As for the ramification filtration on cohomology, the fact that the tame signature has
period $1$ in this case gives that
\begin{eqnarray*}
\lefteqn{\dim_{\Fpbar} \fil^s(H^1(G_K,\Fpbar(\chi))) = }\\
&&\qquad\qquad\quad \left\{ \begin{array}{ll}
0,&\text{if}\,\,\,s < 0,\\
\delta_\triv,&\text{if}\,\,\, 0 \le s < 1 + n_0/(p^2-1),\\
\delta_\triv  + 2 ,&\text{if}\,\,\, 1 + n_0/(p^2-1) \le s < 1 + p/(p-1),\\
\delta_\triv + 2 + \delta_\cyc ,&\text{if}\,\,\, 1 + p/(p-1)) \le s,
\end{array}
\right.
\end{eqnarray*} %
where $\delta_\triv$ (resp.~$\delta_\cyc$) is $1$ or $0$ according to
whether or not $\chi$ is trivial (resp.~cyclotomic).  Unlike Case I,
neither of the spaces $L_i$ can be described in terms of the ramification
filtration, nor can we necessarily detect whether $c_\rho \in L_i$
from the splitting field of $\rho$.

\subsection{Case III}  
Finally we consider Case III, where the tame signature of $\chi$ has the form   %
$(a_0,p)$.   The elements of $W'(\chi,1)$ are then
given in the table:
\begin{equation*}
\begin{tabular}{   >{$}c<{$}|   >{$}c<{$}|   >{$}c<{$}   >{$}c<{$}   >{$}c<{$}   >{$}c<{$}   }\toprule
&J                   &                T            &  \{0\}               &    \{1\}                 & \emptyset  \\ \midrule
\multirow{2}{*}{$1 \le a_0 < p-1$}
&\vec{d}         &           (0,0)          & (p-2,p-1)    & (a_0,p-1)      & (a_0,p) \\
&\vec{b}         &   (a_0,p)         & (a_0+2,p)&(p-1-a_0,1)&(b_0,b_1) \\ \midrule
\multirow{2}{*}{$a_0 = p-1$}
&\vec{d}         &           (0,0)          & (p-1,0)    & (p-1,p-2)      & (0,1) \\
&\vec{b}         &   (p-1,p)               & (1,p-1)    &(p,2)              &(b_0,b_1)\\
\bottomrule
\end{tabular}
\end{equation*}
where 
$$(b_0,b_1) = \left\{ \begin{array}{ll}
(p-2-a_0,p-1), & \mbox{if $a_0 < p-2$,}\\
(p,p-1),&\mbox{if $a_0 = p-2$ or $p=2$,}\\
(p-1,p-2),& \mbox{if $a_0 = p-1$ and $p>2$.}\\
\end{array}\right.$$

As in Case I, there is a unique $J \in S_V(\chi,1)$ for each $V$ in the table, so that $J = J_\max$.
However only $T$, $\{0\}$ and $\emptyset$ are admissible, and the functions
$\delta$ and $\mu$ introduced in Section~\ref{subsec:subspaces} are non-trivial.
Indeed we find that
$\mu(T) = T$, $\mu(\{0\}) = \mu(\{1\})= \{0\}$ and $\mu(\emptyset) = \emptyset$.
If $J=T$, then $L_V^\AH$ is the whole two-dimensional
space $H^1(G_K,\Fpbar(\chi))$, and if $J = \emptyset$, then $L_V^\AH = 0$,
but if $J= \{i\}$ for $i = 0$ or $1$, then $L_V^\AH$ is the same one-dimensional
subspace of $H^1(G_K,\Fpbar(\chi))$ which we will simply denote $L_0$.
We therefore have three possibilities for $W^\AH(\rho)$.

If $1 \le a_0 < p-1$, then $W^\AH(\rho)$ is:
\begin{itemize}
\setlength{\itemindent}{0.2in}
\item[IIIa)] $\{\,V_{(0,0),(a_0,p)}\,\}$ if $c_\rho \not\in L_0$;
\item[IIIb${}_1$)] $\{\,V_{(0,0),(a_0,p)}, V_{(p-2,p-1),(a_0+2,p)}, V_{(a_0,p-1),(p-1-a_0,1)}\,\}$ if $c_\rho \in  L_0$,
         $c_\rho\neq 0$;
\item[IIIc)] $\{\,V_{(0,0),(a_0,p)}, V_{(p-2,p-1),(a_0+2,p)}, V_{(a_0,p-1),(p-1-a_0,1)},
          V_{(a_0,p),(b_0,b_1)}\,\}$ if $c_\rho = 0$.
\end{itemize}

If $a_0 = p-1$, then $W^\AH(\rho)$ is:
\begin{itemize}
\setlength{\itemindent}{0.2in}
\item[IIIa$'$)] $\{\,V_{(0,0),(a_0,p)}\,\}$ if $c_\rho \not\in L_0$;
\item[IIIb${}_1'$)] $\{\,V_{(0,0),(a_0,p)}, V_{(p-1,0),(1,p-1)}, V_{(p-1,p-2),(p,2)}\,\}$ if $c_\rho \in  L_0$,
         $c_\rho\neq 0$;
\item[IIIc$'$)] $\{\,V_{(0,0),(a_0,p)}, V_{(p-1,0),(1,p-1)}, V_{(p-1,p-2),(p,2)},
          V_{(0,1),(b_0,b_1)}\,\}$ if $c_\rho = 0$.
\end{itemize}

Turning to the subspace $L_0$, we now have $n_0 = a_0 + p^2$
and $n_1 = (a_0+1)p$, so that again $n_1' = en_0/(p^2-1)$,
but now $n_0' = n_1' - e$ if $a_0 < p-1$ and $n_0' = e/(p^2-1) = 1$
if $a_0 = p-1$.  Therefore $L_0$ is spanned by the class $c_0$
determined by the formula
$$c_0(E([a]\pi^{n_1'})) = 0,\qquad \displaystyle c_0(E([a]\pi^{n_0'}) ) = 
   \sum_{g\in\gal(l/k)}  \mu^{-1}(g) \tilde{\tau}_0(ga^p)$$
for $a \in l$, where $\tilde{\tau}_0$ is any choice of embedding
extending $\tau_0$.
We thus obtain the criterion that $c_\rho \in L_0$ if and only
if $E([a]\pi^{en_0/(p^2-1)}) \in \ker(c_\rho)$ for all $a \in l$.
In terms of the splitting field $N$ of $\rho$ over $M$, we have
\begin{equation}\label{criterion:norm3}
c_\rho \in L_0 \quad\mbox{if and only if}\quad
E([a]\pi^{en_0/(p^2-1)}) \in \Nm_{N/M}(N^\times)
\mbox{\ for all $a\in l$}.
\end{equation}

As for the ramification filtration, we now have
$$\dim_{\Fpbar} \fil^s(H^1(G_K,\Fpbar(\chi))) = \left\{\begin{array}{ll}
0,&\mbox{if $s < 1 + m/(p^2-1)$,}\\
1,&\mbox{if $1 + m/(p^2-1) \le s < 1 + n_0/(p^2-1)$,}\\
2,&\mbox{if $1 + n_0/(p^2-1) \le s$,}\end{array}\right.$$
where $m = a_0 + 1$ if $1 \le a_0 < p-1$ and $m=1$ if $a_0 = p-1$.
We thus see that 
\begin{equation}\label{criterion:slope3}
L_0 =  \fil^{1+m/(p^2-1)}(H^1(G_K,\Fpbar(\chi)))
 = \fil^{< 1+n_0/(p^2-1)}(H^1(G_K,\Fpbar(\chi))),
 \end{equation}
so that $c_\rho \in L_0$ if and only if $G_K^{n_0/(p^2-1)} \subset \ker(\rho)$.
Moreover since $m < p+1$ and $n_0 > p^2-1$, we have
$L_0   = H^1_\mr(G_K,\Fpbar(\chi)) = H^1_\cg(G_K, \Fpbar(\chi))$,
so that these are precisely the gently ramified classes, which in this case
coincide with the cogently ramified classes.

We remark that if $a_0 = p-2$, we have $L_{V_{(p-2,p-1),(1,1)}}  \subset H^1_\pr(G_K,\Fpbar(\chi))$
by \cite{fontaine}; together with the equality $L_{V_{(p-2,p-1),(1,1)}} = L_{V_{(p-2,p-1),(p,p)}}$ 
provided by \cite{cd}, it follows in this particular case that $W^\AH(\rho) = W(\rho)$.

\subsection{Two irreducible cases}
For completeness, we also list the possibilities when $\rho$ is irreducible.  Recall that in this case we let $W^\AH(\rho) =W(\rho)$
as defined in (\ref{eqn:irred}).

We let $K'$ denote the unramified quadratic extension of $k$ and $k'$
its residue field.  Choose an embedding $\tau':k' \to \Fpbar$ and let
$\psi = \omega_{\tau'}: I_K \to \Fpbar^\times$ denote the associated
fundamental character, so $\psi$ has order $p^4 -1$.  We then have
$$\rho|_{I_K} \sim \begin{pmatrix}\psi^a&0\\0&\psi^{p^2a}\end{pmatrix}$$
for some $a$ with $1 \le a \le p^4 - 1$ and $a \not\equiv 0 \bmod
p^2+1$.  Twisting by characters of $G_K$, we may alter $a$ by multiples
of $p^2 +1$ and hence assume $1 \le a \le p^2$.    Altering our choice
of $\tau'$, we may further assume $a = a_0 + a_1p$ where either
\begin{itemize}
\item[IV)]  $2 \le a_0 \le p -1$ and $1 \le a_1 \le p-2$, or
\item[V)] $1 \le a_0 \le p -1$ and $a_1 = 0$.
\end{itemize}
In Case IV, which is equivalent to $a \not\equiv ip^j \bmod p^2 + 1$
for $i = 1,\ldots,p-1$, $j=0,1,2,3$, we find that
\begin{eqnarray*}
W(\rho) &=&\left\{ \,V_{(0,0),(a_0,a_1)}, V_{(a_0-1,a_1),(p+1-a_0,p-1-a_1)}, V_{(a_0-1,p-1),(p-a_0,a_1+1)},\right.\\
&&\qquad\qquad\qquad\qquad \left. V_{(0,a_1),(a_0-1,p-a_1)}\,\right \}
\end{eqnarray*} %
where the indices in $T$ are ordered so the first embedding is the
restriction of our chosen $\tau'$.
In case V, we find that
\begin{eqnarray*}
W(\rho) &=& \left\{ \, V_{(p-2,p-1),(a_0+1,p)}, V_{(a_0-1,0),(p+1-a_0,p-1)}, V_{(a_0-1,p-1),(p-a_0,1)},\right.\\
&&\qquad\qquad\qquad\qquad \left.V_{(0,0),(a_0-1,p)}\, \right\}
\end{eqnarray*}
with the last weight omitted if $a_0 =1$.

\section{Examples of Galois representations} \label{sec:examples}
We now illustrate the possible behavior discussed in the preceding
section with eight explicit examples for $p=3$, $f=2$.
In the next section, we will exhibit in Table~\ref{8matches} 
numerically matching automorphic
data for each of the Galois representations described here.  
We refer to \cite{ddr} for an extensive collection of examples for more general
$p$ and $f$ and elaboration on methods for obtaining and analyzing them.

We are restricting here to $p=3$, as this is the smallest prime
 for which all the reducible Cases I, II and III arise.
We organize the examples
according to the classification in the preceding section, and we content ourselves
with examples  %
for each type labeled a) or b${}_i$) as these already
illustrate the main new phenomena involving wild ramification in the quadratic case.

In the first two examples, $F = \Q(\sqrt{2})$ while in the last six, $F = \Q(\sqrt{5})$. 
All our representations
$\rho$ take values in $\GL_2(k)$, where $k = \CO_F/3\CO_F$ is viewed as a
subfield of $\overline{\F}_3$ via the embedding labelled $\tau_0$.
We let $\gen$ denote a root of $x^2+2x-1$ if $F = \Q(\sqrt{2})$, and a root of
$x^2 - x -1$ if $F = \Q(\sqrt{5})$, and in either case we use the same symbol
$\bargen$ for its image in $k \subset \overline{\F}_3$.  

In the list below,
 we describe $\rho$ by specifying its projective splitting field, its conductor  %
(prime to $3$) and its local behavior at $p=3$ up to an unramified quadratic twist.
In each of our examples one can show there is a
unique representation $\rho$ satisfying this description, except for those in Case III, where there
are two such representations differing by a quadratic twist. %

\subsection{Case I} We use examples with tame signature $(a_0,a_1) = (1,2)$,
so $n_0 = 7$ and $n_1 = 5$.  This means that $\rho|_{G_K}$ is a twist of a 
representation of the form 
$$\begin{pmatrix}\chi&*\\0&1\end{pmatrix},  \qquad
\mbox{with $\chi|_{I_K} = \omega_{\tau_0}^7 = \omega_{\tau_1}^5$;}$$
in particular $\chi$ is primitive and generic. 

In all our examples,
 $\chi$ will in fact have the form $\omega_\pi^7: \gal(M/K) \to k^\times$
where $\pi^8$ is a uniformizer of $K$, $M= K(\pi)$ and $\omega_\pi$ is the associated
fundamental character.  The class $c_\rho \in H^1(G_K,k(\omega_\pi^7))$ thus
corresponds via local class field theory to a $\gal(M/K)$-linear homomorphism
$$M^\times  \longrightarrow \gal(N/M) \cong k(\omega_\pi^7)$$
with kernel $\Nm_{N/M}(N^\times)$ where $N$ is the projective splitting field of
$\rho|_{G_K}$.   This kernel contains $IM^\times$,
where $I$ is the kernel of the surjection $\Z[\gal(M/K)] \to k$ induced
by $\omega_\pi^7$.  As a $k$-vector space $M^\times/IM^\times$
is two-dimensional, spanned by $E_3(\pi^5) \equiv 1 + \pi^5$ and
$E_3(\pi^7) \equiv 1 + \pi^7$, and $\Nm_{N/M}(N^\times)/IM^\times$
is a one-dimensional subspace that determines $W^\AH(\rho|_{G_K})$
via (\ref{criterion:norm}).

\subsection*{Example Ia}
Let $F = \Q(\sqrt{2})$ and let $E$ denote the
splitting field over $\Q$ of the polynomial
$$ f_{\rm Ia}(x) = x^{10}-24 x^7-42 x^6-24 x^4-48 x^3-18 x^2-32 x-96.$$
Then $F \subset E$, and there is an isomorphism $\varrho:\gal(E/F) \to \PGL_2(k)$
that lifts to a representation $\rho:G_F \to \GL_2(k)$ of conductor $\gp_2^6$,
where $\gp_2 = (\sqrt{2})$.

Up to an unramified quadratic twist, 
the local representation $\rho|_{G_K}$ has the form 
$$\omega_\pi^5 \otimes \begin{pmatrix}\omega_\pi^7&*\\0&1\end{pmatrix}$$
where $\pi^8 = 3$, and the splitting field $N$ of the
projective local representation is that of the polynomial
$x^9 + 6x^7 + 3x^6 + 6.$  
Here and in the later examples, we are using the database 
described in \cite{jr} to pass from the global polynomial
to a local $3$-adic Eisenstein polynomial.  {As will be explained in more detail in
\cite{ddr}, the above form for the local representation is determined up to twist by the maximal
tamely ramified subfield of $N$ (in this case $M = K(\pi)$), the action of $\gal(M/K)$
on $\gal(N/M)$ and the choice of isomorphism $\varrho$.  The twist is then specified, up to an
unramified quadratic character, as part of the data characterizing the lift $\rho$ of
the projective representation $\varrho$.  Using a  {\em Magma} program described
in \cite{ddr},} we find that $\Nm_{N/M}(N^\times)/IM^\times$
consists of the classes of
elements of the form $1 + [a]\pi^5 - [a]^3\pi^7$ for $a \in k$.
Taking into account the twist by $\omega_\pi^5$, we conclude from
(\ref{criterion:norm}) that $W^\AH(\rho|_{G_K}) = \{ V_{(2,1),(1,2)} \}$.

\subsection*{Example Ib${}_1$}
 Let $F = \Q(\sqrt{2})$ again and let $E$ denote the
splitting field over $F$ of the polynomial
$$f_{\rm Ib_1}(x) =x^{10}-9 x^8+78 x^6-246 x^4-48 x^3+459 x^2+224 x-75.$$
We again have $F \subset E$ and an isomorphism $\gal(E/F) \cong \PGL_2(k)$
lifting to a representation $\rho:G_F \to \GL_2(k)$ of conductor $\gp_2^6$.

Up to an unramified quadratic twist, $\rho|_{G_K}$ has the form 
$$\omega_\pi^2 \otimes \begin{pmatrix}\omega_\pi^7&*\\0&1\end{pmatrix}$$
where now $\pi^8 = - 3$ and
the local projective splitting field $N$ is that of
$x^9 +3x^7 + 3.$
In this case however $\Nm_{N/M}(N^\times)/IM^\times$
consists of the classes of $1 + [a]\pi^5$ for $a \in k$, so
(\ref{criterion:norm}) implies that $c_\rho \in L_1$.  Taking into account
the twist by $\omega_\pi^2$, we conclude that
$W^\AH(\rho|_{G_K}) = \{ V_{(2,0),(1,2)}, V_{(0,0),(2,3)} \}$.

\subsection*{Example Ib${}_2$} Now, 
and for all the remaining examples, let $F =\Q(\sqrt{5})$.  Let $E$ denote the
splitting field over $\Q$ of the polynomial $f_{\rm Ib_2}(x) =$
 $$x^{10}-2 x^9+9 x^8+48 x^7-132 x^6+504 x^5+228 x^4-1824 x^3 +6894 x^2-7676 x+4462.$$
We again have $F \subset E$ and an isomorphism $\gal(E/F) \cong \PGL_2(k)$
lifting to a representation $\rho:G_F \to \GL_2(k)$ of conductor $(2)^5$.

Up to an unramified quadratic twist, $\rho|_{G_K}$ has the form 
$$\omega_\pi^7 \otimes \begin{pmatrix}\omega_\pi^7&*\\0&1\end{pmatrix}$$
where $\pi^8 = - 3$ and
the local projective splitting field $N$ is that of
$x^9 + 6x^5 + 6.$
We now find that $\Nm_{N/M}(N^\times)/IM^\times$
consists of the classes of $1 + [a]\pi^7$ for $a \in k$, so that
 $c_\rho \in L_0$ by (\ref{criterion:norm}) 
(or by (\ref{criterion:slope1}) since $G_K^{7/8} \subset \ker(\rho)$).
 Taking into account the twist by $\omega_\pi^7$, we conclude that
 $W^\AH(\rho|_{G_K}) = \{ V_{(1,2),(1,2)}, V_{(1,1),(2,1)} \}$.

\subsection{Case II} We use examples with tame signature $(a_0,a_1) = (1,1)$,
so $n_0 =  n_1 = 4$.  Thus $\rho|_{G_K}$ is a twist of a 
representation of the form 
$$\begin{pmatrix}\chi&*\\0&1\end{pmatrix},  \qquad
\mbox{with $\chi|_{I_K} = \omega_{\tau_0}^4 = \omega_{\tau_1}^4$;}$$
in particular $\chi$ is imprimitive and generic.

Note that we may write $\chi = \mu\omega_\pi$ where $\mu$ is unramified
and $\pi^2 = -3$, so $\omega_\pi$ is the cyclotomic character.  In all our
examples, we will have $c_\rho \in H^1_\ar(G_K,k(\chi))$, and since
$H^1_\ur(G_K,k(\chi)) = 0$, we see that $W^\AH(\rho|_{G_K})$ is determined
by whether $c_\rho$ is a multiple of either of the classes defined in
(\ref{criterion:alignment}).

\subsection*{Example IIa}  Let $E$ be the splitting field over $\Q$
of the polynomial
$$f_{\rm IIa}(x) = x^4-x^3+2 x-11.$$
A representation $\rho:G_F \to \GL_2(\F_3)$ of conductor $(7)$
with projective splitting field $E$ is given by the $3$-torsion of a quadratic twist of
the base-change to $F = \Q(\sqrt{5})$ of the elliptic curve over $\Q$ with Cremona label 175A.

Up to an unramified quadratic twist,
the local representation $\rho|_{G_K}$ has the form 
$$ \begin{pmatrix}\omega_\pi&*\\0&1\end{pmatrix}$$
where $\pi^2 = - 3$, so we can take $M = K(\pi)$ and $l = k$.
The local projective splitting field $N$ is that of
$x^3 +3x + 3.$
Since $\omega_\pi$ is cyclotomic, we have that 
$H^1_{\pr}(G_K,\F_3(\omega_\pi)) = H^1_{\ar}(G_K,\F_3(\omega_\pi))$ has codimension one
in $H^1(G_K,\F_3(\omega_\pi))$.
One can check directly that $G_K^{3/2} \subset \ker(\rho)$ and hence that
$c_\rho \in H^1_{\ar}(G_K, \F_3(\omega_\pi))$, or deduce this from the fact
$\rho$ is defined by an elliptic curve with good ordinary reduction at $3$.
On the other hand since $c_\rho$ is non-trivial and takes values in $\F_3$,
but the homomorphism in (\ref{criterion:alignment}) is simply $\tau_i$
and hence has image of order $9$, it follows that $c_\rho \not\in L_0 \cup L_1$.
Therefore $W^\AH(\rho|_{G_K}) = \{ V_{(0,0),(3,3)}, V_{(0,0),(1,1)} \}$.

\subsection*{Example IIb${}_1$}
 Let  $E$ be the splitting field over $F = \Q(\sqrt{5})$ of %
$$f_{{\rm IIb}_1}(x) = x^6 - 3\gen x^5 + 3\gen x^4 + (6\gen + 6)x^3 - (21\gen +12)x^2 + (21\gen + 12)x - 8\gen - 4.$$%
The Shimura curve associated to the units of a maximal order of a
quaternion algebra over $F$ ramified at one archimedean place and 
the prime $\gp_{61} = (3-7\gen)$ has genus two, and its
Jacobian has real multiplication by $F$ (see \cite[Remark~3]{lassina:expmath}).
The $3$-torsion points of this Jacobian give rise to a representation $\rho:G_F \to \GL_2(k)$
of conductor $\gp_{61}$ with $E$ as its projective splitting field.  Note that unlike the preceding examples, $E$ is not
Galois over $\Q$.

Up to an unramified quadratic twist, the local representation  $\rho|_{G_K}$ has the form 
$$ \nu^{-1}\otimes \begin{pmatrix}\nu^2\omega_\pi&*\\0&1\end{pmatrix}$$
where $\pi^2 = - 3$ and $\nu$ is the unramified character of $G_K$ sending
$\frob_K$ to (the reduction of) $\gen^3$.  
We let $M = L(\pi)$ where $L$ is the unramified extension of $K$ of degree~$4$,
so also $[l:k] = 4$.  The splitting field of the character $\nu^2\omega_\pi$ is not
of the form required for the construction of Section~\ref{subsec:basis}, so we have
adjoined $\pi$ in order to obtain a field of the required form; note that the extension
$M/K$ is not cyclic.  Note also that since $\nu^2$ is non-trivial, we have
$H_\ar^1(G_K, k(\nu^2\omega_\pi)) =  H^1(G_K, k(\nu^2\omega_\pi))$.
The class $c_\rho$ now corresponds to a $\gal(M/K)$-linear homomorphism
$$M^\times/IM^\times  \longrightarrow \gal(N/M) \cong k(\nu^2\omega_\pi)$$
where $I$ is the kernel of the surjection $\Z[\gal(M/K)] \to k$ induced
by $\nu^2\omega_\pi$, and $N$ is the composite of $M$ with the projective
local splitting field of $\rho$.   As a $k$-vector space, $M^\times/IM^\times$ is
two-dimensional, consisting of the classes of $E_3([a]) \equiv 1 + [a]\pi$
for $a$ in the kernel of $\tr_{l/k'}$, where $k'$ is the quadratic extension of $k$.
Unravelling (\ref{criterion:alignment}), we find that $c_\rho \in L_0$ (resp.~$L_1$)
if and only if $c_\rho$ is trivial on those $1 + [a]\pi$ such that $a^8 = \nu^2(\frob_K) = \gen^2$
(resp.~$a^8 = \nu^6(\frob_K) = -\gen^2$).  Explicit computation of elements of
$\Nm_{N/M}(N^\times)$ shows that indeed $c_\rho$ is in $L_1$ (and hence
not in $L_0$ since $c_\rho \neq 0$), so that 
$W^\AH(\rho|_{G_K}) = \{ V_{(0,0),(3,3)}, V_{(0,0),(1,1)}, V_{(0,2),(2,2)} \}$.

\subsection*{Example IIb${}_2$} We may take Example IIb${}_1$ and replace $\rho$
by $\rho\circ \sigma$, where $\sigma$ is the outer automorphism of $G_F$ induced by
conjugation by an element of $G_\Q$ extending the non-trivial element of $\gal(F/\Q)$.
The resulting representation has conductor $\gp_{61}' =
(4-7\gen)$ and projective splitting field $\sigma(E)$;
the character $\nu$ in the description of $\rho|_{G_K}$ is the same
as in Example~IIb${}_1$, but the kernel of the homomorphism induced
by $c_\rho$ would be replaced by its Galois conjugate.  We therefore
conclude that $c_\rho$ is in $L_0$ instead of $L_1$, so that
$W_\AH(\rho|_{G_K}) = \{ V_{(0,0),(3,3)}, V_{(0,0),(1,1)}, V_{(2,0),(2,2)} \}$.
The corresponding system of
Hecke eigenvalues is obtained from the one in Example IIb${}_1$ by
interchanging each $a_v$ with $a_{\sigma(v)}$. 
(Note that a similar procedure could not have been used to generate an example
of type Ib${}_2$ from Ib${}_1$ since the inequality $n_1 < n_0$ would not
be preserved.)

Alternatively, we could obtain an example of type IIb${}_2$ by replacing the
representation $\rho$ in Example IIb${}_1$ by its composite with the automorphism of $\GL_2(k)$ induced by $\frob$ on $k$.  The projective splitting field is then the same as in
Example IIb${}_1$, as is the description of $\rho|_{G_K}$, except that $\nu$ is
replaced by the unramified character sending $\frob_K$ to $\alpha$ and the
homomorphism corresponding to $c_\rho: M^\times \to k(\nu^2\omega_\pi)$
is obtained from the preceding one by composing with $\frob$.  Note that $N$ and
$IM^\times$ do not change, but the criteria for $c_\rho$ to be in $L_0$ and $L_1$
in terms of $\Nm_{N/M}(N^\times)$ are interchanged.  In this case the corresponding
system of Hecke eigenvalues is obtained from the one
in Example IIb${}_1$ by replacing each $a_v$ with $\frob(a_v)$. 
Finally of course, we could just as well have obtained an
example of type IIb${}_1$ by replacing the original $\rho$ with
$\frob\circ\rho\circ\sigma$.

\subsection{Case III} We use examples with tame signature $(a_0,a_1) = (1,3)$,
so $n_0 = 10$ and $n_1 = 6$.   Thus $\rho|_{G_K}$ is a twist of a 
representation of the form 
$$\begin{pmatrix}\chi&*\\0&1\end{pmatrix},  \qquad
\mbox{with $\chi|_{I_K} = \omega_{\tau_0}^2 = \omega_{\tau_1}^6$;}$$
in particular $\chi$ is primitive and non-generic. 

In both our examples $\chi$ will in fact have the form $\omega_\pi: \gal(M/K) \to k^\times$
where $\pi^4$ is a uniformizer of $K$, $M= K(\pi)$ and $\omega_\pi$ is the associated
fundamental character.  The class $c_\rho$ will be non-trivial, so that
$W^\AH(\rho|_{G_K})$ is determined by whether $c_\rho$ lies in the
space $L_0$ described in (\ref{criterion:norm3}) or (\ref{criterion:slope3}).
Note also that since $a_0 = p-2$, we know in fact in this case that
 $W^\AH(\rho|_{G_K}) = W(\rho|_{G_K})$ by the remark at the end
 of Section~\ref{sec:quadratic}.

For both examples, there are in fact two representations with the
given description; choosing either to be $\rho$, the other is $\delta\otimes\rho$
where $\delta$ the non-trivial character of $\gal(F(\zeta_5)/F)$.
\subsection*{Example IIIa}

Let $E$ denote the
splitting field over $\Q$ of the polynomial
 $$f_{\rm IIIa}(x) = x^{10}-5 x^9+135 x^6-360 x^5+405 x^4-270 x^3+135 x^2-45 x+9.$$
We then have $F = \Q(\sqrt{5}) \subset E$ and an isomorphism %
$\gal(E/F) \cong \PSL_2(k)$ lifting to a representation
$\rho:G_F \to \GL_2(k)$ of conductor $\gp_5^3$, where $\gp_5 = (\sqrt{5})$.

Up to an unramified quadratic twist,
$\rho|_{G_K}$ has the form 
$$\omega_\varpi \otimes \begin{pmatrix}\omega_\varpi^2&*\\0&1\end{pmatrix}$$
where $\varpi^8 = 3\omega^2$, so we may take $M = K(\pi)$ 
with $\pi = \varpi^2$.   The splitting field $N$ of the
projective local representation is that of the polynomial
$x^9 + 9x + 6,$
and we find that $G_K^{5/4} \not\subset \ker(\rho)$, so that $c_\rho \not\in L_0
 = \fil^{5/4} H^1(G_K,\overline{\F}_3(\omega_\pi))$ by (\ref{criterion:slope3}).
 Taking into account the twist by $\omega_\varpi$, we conclude that
 $W^\AH(\rho|_{G_K}) = \{ V_{(1,0),(1,3)} \}$.

\subsection*{Example IIIb${}_1$} 
Let $E$ denote the
splitting field over $\Q$ of the polynomial
 $$f_{\rm IIIb_1}(x) =x^6-3 x^5+5 x^3-5.$$
We again have $F \subset E$ %
and an isomorphism $\gal(E/F) \cong \PSL_2(k)$ lifting to a representation
$\rho:G_F \to \GL_2(k)$, now of conductor $(2)\gp_5^3$.

Up to an unramified quadratic twist, $\rho|_{G_K}$ again has the form 
$$\omega_\varpi \otimes \begin{pmatrix}\omega_\varpi^2&*\\0&1\end{pmatrix},$$
but now $\varpi^8 = 3$, and we take $M = K(\pi)$ with $\pi = \varpi^2$.   
In contrast to the preceding example, we find that $G_K^{5/4} \subset \ker(\rho)$,
so that $c_\rho \in L_0$ by (\ref{criterion:slope3}).
Taking into account the twist by $\omega_\varpi$, we conclude that
 $W^\AH(\rho|_{G_K}) = \{ V_{(1,0),(1,3)}, V_{(0,0),(3,3)}, V_{(0,0),(1,1)} \}$.

\section{Numerical matching with automorphic forms}     \label{sec:numerical}
To facilitate computations, both here and in the sequel \cite{ddr}, we work
with algebraic automorphic forms on definite 
quaternion algebras over totally real fields.  Recall that these are related to Hilbert modular forms
by the Jacquet--Langlands correspondence, and under mild hypotheses
give the set of weights of forms giving rise to $\rho$ in the sense of 
Conjectures~\ref{conj:serre1} and~\ref{conj:serre2}.

More precisely, we will consider totally real fields $F$ in which $p$ is inert and
definite quaternion algebras $B$ over $F$ which are split at $p$, and
we present pairs $(\phi,\rho)$ where:
\begin{itemize}
\item $\phi =  (a_v, d_v)_{v\in \Sigma_\phi}$ is a system of eigenvalues for the standard Hecke operators $T_v$ and $S_v$ (as defined in \cite{rlt:inv}) acting on mod $p$ algebraic
modular forms for $B$ of some level $\gn_\phi$ (where $\Sigma_\phi$ is a large set of good primes);
\item $\rho:G_F \to \GL_2(\Fpbar)$ is a Galois representation unramified outside $p\gn_\phi$
such that $\rho(\frob_v)$ has characteristic polynomial $x^2  \pm a_v x  + d_v \N(v)$ for all
$v \in \Sigma_\phi$;
\item the set of weights for which $\phi$ occurs at level $\gn_\phi$ is precisely
$W^\AH(\rho|_{G_{F_p}})$.
\end{itemize}
The reason for the sign ambiguity in the trace of $\rho(\frob_v)$ is that in practice
we work with the associated projective representation.
The $\rho$ we consider, particularly in \cite{ddr}, %
are typically constructed independently from automorphic forms.  The existence
of a numerically matching $\phi$ can be viewed as evidence for the modularity part of 
Conjecture~\ref{conj:serre2} (and hence Conjecture~\ref{conj:serre1}).

    For each of the eight Galois representations $\rho$ from the previous
section, we exhibit a corresponding $\phi$ here, taking $\Sigma_\phi$ 
to be all good primes with norm at most $100$.  The methods for computing $\phi$
are based on those described in~ \cite{lassina:mathcomp} and Appendix~B
of~\cite{breuil:appendix}.

 \subsection{A summarizing table}   Table~\ref{summarizing}
 summarizes our eight examples, adding some more information.  
 Note that in all cases
  besides the conjugate cases IIb${}_1$ and IIb${}_2$, the polynomial 
  $F_c(x) := f_c(x)$ has coefficients in $\Q$.    Its Galois group is given in the $G$ column.  %
 In this column, an exceptional
 isomorphism identifies the group $\PSL_2(9)$ with the alternating group $A_6$.  The group $\mathrm{P\Gamma L}_2(9) = \Aut(\PSL_2(9))$ contains
 $\PSL_2(9) \cong A_6$ with index four and the three intermediate groups are $M_{10}$, $\PGL_2(9)$, and $S_6$.  
 The entry shared by the IIb${}_1$ and IIb${}_2$ rows is the Galois group
  of the product $F_{{\rm IIb}}(x) := f_{{\rm IIb}_1}(x) f_{{\rm IIb}_2}(x) \in \Z[x]$.  The $D$ column gives
  the field discriminant of $\Q[x]/F_c(x)$.     The largest slope $s$ is 
  explained in the next subsection. 
    
  \begin{table}[htb]\small
\setlength{\tabcolsep}{2pt}
\begin{tabular}{  >{$}c<{$}  >{$}c<{$}|  >{$}c<{$}  >{$}c<{$}  >{$}c<{$}|  >{$}c<{$}  >{$}r<{$}  >{$}c<{$}|  >{$}l<{$} }\toprule
 \text{Ex.} & F   &  G &  D &   s &  \mathfrak{n}  & \N(\mathfrak{n})  & \mathfrak{f} &\multicolumn{1}{c}{\text{Weights}\,\, ($W^{\AH}(\rho)$)} \\\midrule
\mathrm{Ia} & \Q(\sqrt{2}) & \mathrm{P \Gamma L}_2(9)    &  2^{27} 3^{15} &15/8  &  \gp_2^6  & 64 &  \gp_2^6  &  [2,1;1,2]    \\

\mathrm{Ib}_1 &  \Q(\sqrt{2})  &  \mathrm{P \Gamma L}_2(9)    & 2^{28} 3^{15} &15/8  &\mathfrak{p}_2^6   & 64  &  \mathfrak{p}_2^6  & [2,0;1,2], \mathbf{[0,0;2,3]}  \\

\mathrm{Ib}_2 & \Q(\sqrt{5})  &  \mathrm{P \Gamma L}_2(9)    &  2^{31} 3^{13} 5^3  &13/8& (2)^5   & 1024  &(2)^5& [1,2;1,2],  \mathit{[1,1;2,1]}  \\\midrule

\mathrm{IIa}  & \Q(\sqrt{5})    &  S_4 & -3^3 5^2 7^2      &3/2   & (7) & 49    & \mathcal{O}    &  [0,0;3,3], \mathit{[0,0; 1,1]}  \\

\mathrm{IIb}_1 & \Q(\sqrt{5})  &  \multirow{2}{*}{$A_6^2.2$}     & \multirow{2}{*}{$3^{12} 5^6 61^2$} & \multirow{2}{*}{3/2}  & \mathfrak{p}_{61} & 61 &\mathcal{O}&  [0,0;3,3],
\mathit{[0,0; 1,1]}, \mathbf{[0,2; 2,2]}  \\

\mathrm{IIb}_2 & \Q(\sqrt{5})  &                 &      &  & \mathfrak{p}_{61}' & 61 &\mathcal{O}&  [0,0;3,3], \mathit{[0,0; 1,1]},  \mathbf{[2,0; 2,2]}  \\\midrule

\mathrm{IIIa}   &  \Q(\sqrt{5})   &  M_{10}   & 3^{18} 5^{10}    &9/4& \mathfrak{p}^3_5 &  125 & \mathcal{O} & [1,0;1,3]  \\

\mathrm{IIIb}_1 &  \Q(\sqrt{5})  &   S_6      & 2^2 3^6 5^2  &5/4&  (2) \mathfrak{p}^3_5 & 500  &  \mathcal{O}  & [1,0;1,3] ,  \mathit{[0,0; 3,3]},  \mathit{[0,0;1,1]} \\
\bottomrule
\end{tabular}
\vspace{5pt}
\caption{\label{summarizing}Information on the eight examples.  The weight in ordinary type 
is computed from the tame signature.  As explained in \S\ref{slopes}, weights
in italics come from small slopes and weights in boldface come
from other sources.}  
\end{table}

\subsection{Slopes} \label{slopes}  For a separable polynomial $f(x) \in \Q_p[x]$, wild ramification
in the algebra $A=\Q_p[x]/f(x)$ can be measured by slopes, as explained in \cite{jr}.  These slopes are
breaks in the upper numbering of {\cite[IV.3]{serre_cl}}, increased by $1$.   When all factors of $f(x)$ have degree
$\leq 11$, they are computed automatically by the website of \cite{jr}.

A common situation in our current setting is that $f(x) \in \Q[x]$ has degree ten, and factors
over $\Q_3$ into a primitive nonic and a linear factor, giving $A = B \times \Q_3$.    In this case, 
the primitive nonic field $B$ has a certain largest slope $s$ with multiplicity eight and $0$ with multiplicity one.  As
$\Q_3$ has the trivial slope $0$ as well,  $\ord_3(D) = 8s$.   This situation occurs
in our four cases with $F_c(x)$ decic, namely Ia, Ib${}_1$, Ib${}_2$, and IIIa.   
The other cases are similiar.  For example, $\Q_3[x]/f_{{\rm IIIb}_1}(x)$ is a sextic field 
with a tame subfield of degree two.  In this case, $\ord_3(D)=6$ decomposes 
as $4s+1+0$; the $1$ comes from the tame subfield and $s=5/4$ is the 
quantity of current interest.  

The slope column illustrates that some extra weights  
come simply from $s$ being smaller than the maximum allowed
by the tame signature.   For example, for the tame signature $(a_0,a_1)=(2,2)$,
the maximum allowed $s$ is $5/2$, while our examples are peu ramifi\'ee and have
slope $3/2$.      However other extra
weights are not simple consequences of small slopes.
The sequel paper \cite{ddr} will illustrate a principle clear from the theory
here: as the local degree $[K:\Q_p]$ increases, slopes account for a 
decreasing fraction of the 
phenomenon of extra weights.   

\subsection{The class set $\PGL_2(9)^\natural$} 
Table~\ref{matching} summarizes how one does projective matching for general
$\rho$ into $\GL_2(k)$ where $k$ has order $9$.   On the automorphic side, one has the pairs $(a_v,d_v) \in k \times k^\times$.
On the Galois side, the most immediately available quantities are partitions $\lambda_v$
with parts being the degrees of the irreducible factors of $f_c(x)$ in the completed
ring $\cO_v$.  
  
\begin{table}[htb]\small
{\setlength{\tabcolsep}{5pt}
\begin{tabular}{>{$}c<{$}|  >{$}c<{$}  >{$}c<{$}  >{$}c<{$}  >{$}c<{$}  >{$}c<{$}  >{$}c<{$} | >{$}c<{$}  >{$}c<{$}  >{$}c<{$}  >{$}c<{$}  >{$}c<{$}  }
\toprule
d_v & \multicolumn{6}{c|}{\mbox{$d_v$ is a square}} & \multicolumn{5}{c}{\mbox{$d_v$ is not a square}} \\
\midrule 
b_v  & 1 &  0 & 1 &  2 & \w^2 & \w^6 & 0 & \w & \w^3 & \w^5 & \w^7 \\
\midrule
\mathrm{PGL}_2(9)^\natural & 1 & 2u &   3 &  4 &   5A        &  5B  & 2v & 8A & 8B & 10B & 10A \\
\midrule
\mathrm{PGL}_2(9)   & 1^{10} & 2^4 \, 1^2 & 3^3 \, 1 &   4^2  \, 1^2 & 5^2 &   5^2 & 2^5 &  8 \, 1^2 & 8 \, 1^2 & 10 & 10 \\
A_6          &1^6 &  2^2 \, 1^2 & \!\! 3^2 \mbox{ or } 3 \, 1^3 \! &     4 \, 2 &   5 \, 1&    5 \, 1 & &  \\
S_4          &1^4 & \!  2^2 \mbox{ or } 2 \, 1^2 \!\!  & 3 \, 1 &      4   &  &\\
\bottomrule
\end{tabular}
}
\vspace{5pt}
\caption{\label{matching} The class set $\mathrm{PGL}_2(9)^\natural$ and its view from the
automorphic and Galois sides.  The outer involution $\overline{\cdot}$ of $\mathrm{PGL}_2(9)$ makes the interchanges
$5A \leftrightarrow 5B$, $8A \leftrightarrow 8B$, and 
$10A \leftrightarrow 10B$, and fixes the other five classes. }
\end{table}

When $f_c(x)$ is chosen to be a decic in the standard way,
the table explains how the projective quantities $b_v = a_v^2/d_v \bN(v)$
correlate with the decic partitions $\lambda_v$.  
In fact, let $\PGL_2(9)^\natural$ be the set of conjugacy 
classes in the group $\PGL_2(9)$.   Then the Frobenius 
class $\Fr_v \in \PGL_2(9)^\natural$ determines both
$b_v$ and $\lambda_v$.  Conversely,  the pair
$(b_v,\lambda_v)$ determines $\Fr_v$.  
Using the \verb@FrobeniusElement@ command \cite{doks} in {\em Magma},
with adaptations to account for ground field $F$ rather than $\Q$, we have gone
beyond partitions
and have in all cases identified the correct label $A$ or $B$, directly
from the polynomial $f_c(x)$.    Table~\ref{matching} also 
has lines corresponding to our sometimes replacing
decic polynomials by sextic and
quartic polynomials.

\subsection{Matching for our eight examples}
Table~\ref{8matches} is headed by the ten smallest split primes $p$ for each of the two fields $F$ in question.   For each
$p$, it gives one of the two $v$ above it.  The conjugate prime $\sigma(v)$ is obtained by the substitution $\alpha \mapsto -\alpha-2$
in the case $F=\Q(\sqrt{2})$ and $\alpha \mapsto 1-\alpha$ in the case $F = \Q(\sqrt{5})$.
\begin{table}[htb]
\small
\setlength{\tabcolsep}{2.4pt}
\begin{tabular}{  >{$}c<{$}   >{$}c<{$}|   >{$}c<{$}   >{$}c<{$}   >{$}c<{$}   >{$}c<{$}   >{$}c<{$}   >{$}c<{$}   >{$}c<{$}   >{$}c<{$}   >{$}c<{$}   >{$}c<{$}   }
\multicolumn{3}{c}{$F = \Q(\sqrt{2})$}&\multicolumn{9}{c}{} \\ 
\toprule
&p&7 & 17 & 23 & 31 & 41 & 47 & 71 & 73 & 79 & 89  \\
\multicolumn{2}{c|}{$p \mod  3$} & 1 & 2 & 2 & 1 & 2 & 2 & 2 & 1 & 1 & 2 \\
& v & 1 + 2 \w & 2 + 3 \w &  4 - \w &  3 + 4 \w & 5 - 2 \w & 6 - \w & {7 + 6 \w }& 7 - 2 \w &8 - \w & {10+7 \w} \\
\midrule
\multirow{3}{*}{$\mathrm{Ia}$} & a_v& \w^3 & \w^3 & 0 & \w^3 & \w & \w^2 & 2 & \w^3 & 1 & \w^2 \\
                  &d_v &  \w^7 & \w^7 & \w^6 & \w^3 & \w^7 & \w^6 & \w^2 & 1 & \w^3 & \w  \\
& \Fr_v&10A & 8B & 2u & 8B & 10A & 5A & 5A & 5B & 10B & 10A  \\
\midrule
\multirow{3}{*}{$\mathrm{Ib}_1$} &a_v&1 & \w & \w^2 & \w^6 & \w & \w^3 & 0 & \w^2 & \w^7 & 2  \\
                &d_v &\w & \w^5 & \w^6 & \w & \w & \w^2 & \w^2 & 2 & \w & \w^7 \\
&\Fr_v&10A & 8A & 5A & 8B & 10B & 3 & 2u & 3 & 10B& 10B  \\
\bottomrule
\multicolumn{12}{c}{} \\
\multicolumn{12}{c}{} \\
\multicolumn{3}{c}{$F = \Q(\sqrt{5})$}&\multicolumn{9}{c}{} \\
\toprule
&     p &            11 & 19 & 29 & 31 & 41 & 59 & 61 & 71 & 79 & 89 \\
\multicolumn{2}{c|}{$p \mod 3$} & 2 & 1 & 2 & 1 & 2 & 2 & 1 & 2 & 1 & 2 \\
&  v & 2 - 3 \w & 1 - 4 \w & 5 + \w & 3 - 5 \w & 6 + \w & 2 - 7 \w & 4 - 7 \w & 8 + \w &  5 -  8 \w & 10 - \w \\
\midrule
\multirow{3}{*}{$\mathrm{Ib}_2$}&a_v&\w^2 & \w^6 & \w & \w^7 & \w & 2 & \w & \w & \w^7 & \w^2 \\
&d_v& \w^3 & \w^2 & 2 & \w & \w^5 & \w^3 & \w^3 & 2 & \w & \w^7 \\
 &  \Fr_v & 10B & 5A & 5A & 10B & 8A & 8A & 10A & 5A & 10B & 8A \\
\midrule
\multirow{2}{*}{$\mathrm{IIa}$}  &a_v&0 & 0 & 2 & 2 & 2 & 2 & 1 & 1 & 2 & 0 \\
         &  \Fr_v   &     2u & 2u & 4 & 3 & 4 & 4 & 3 & 4 & 3 & 2u \\
 \midrule
\multirow{2}{*}{$\mathrm{IIb}_1$} &a_v& \w^7 & 0 & \w^6 & \w^5 & \w & \w^3 & 0 & 2 & \w^3 & 1 \\
 &\Fr_v&                  5A & 2u & 3 & 5A & 5B & 5A & 2u & 4 & 5B& 4 \\
 \midrule
\multirow{2}{*}{$\mathrm{IIb}_2$} &a_v& \w^6 & \w^6 & \w^3 & \w^2 & \w^7 & \w &  & \w^5 & \w^6 & \w^6 \\
 &\Fr_v &                  3 & 4 & 5A & 4 & 5A & 5B &  & 5B & 4& 3 \\
\midrule
\multirow{2}{*}{$\mathrm{IIIa}$} & a_v  &  0 & \pm\w & \pm2 & \w^5 & \w^5 & \pm\w^6 & \w^2 & \w^5 & 0 & \pm\w^7 \\
 & \Fr_v  & 2u & 5A & 4 & 5A & 5B & 3 & 4 & 5B & 2u & 5A    \\
\midrule
\multirow{2}{*}{$\mathrm{IIIb}_1$}&a_v &2 & \pm\w^5 & \pm{\w^7} & \w^7 & {\w^5} &\pm \w^6 & \w^6 & \w & \pm\w^3 & \pm1 \\
&  \Fr_v   & 4& 5A & {5A} & 5B & {5B} & 3 & 4 & 5B & 5B & 4 \\
\bottomrule
\end{tabular}
\vspace{5pt}
\caption{\label{8matches} Matching in our eight examples for ten $v$} 
\end{table}
For each example, we list the classes $\Fr_v \in \PGL_2(9)^\natural$ associated to $\rho$
and the eigenvalues $a_v$ and $d_v$ of a numerically matching eigenform $\phi$.
We omit the lines for $d_v$ when they are identically $1$.
Recall that in Examples~IIIa and IIIb${}_1$ there are two choices for $\rho$ differing
by twist by the quadratic character $\delta: \gal(\Q(\zeta_5)/F) \to \{\pm 1\}$;
accordingly we list the two matching eigenforms, each obtained from the other
by replacing $a_v$ with $\delta(v)a_v$ where 
$\delta(v) = \delta(\Fr_v) = \left(\frac{\N_{F/\Q}(v)}{5}\right)$.

For each $v$ listed in the table, the eigenvalues $a_{\sigma(v)}$ and $d_{\sigma(v)}$
can be recovered as follows:  In all the examples, one has $d_{\sigma(v)} = d_v^3$, and
in all but IIb${}_1$, IIb${}_2$ and IIIb${}_1$, one has $a_{\sigma(v)} = a_v^3$.
In examples IIb${}_1$ and IIb${}_2$, 
with eigenvalues $a_v'$ and $a_v''$ respectively, one has
 $a'_{\sigma(v) }= a''_v$;
 finally in IIIb${}_1$, one has $a_{\sigma(v)} = \delta(v)a_v^3$.
 Similarly, in all the examples but IIb${}_1$ and IIb${}_2$, one has 
 $\Fr_{\sigma(v) }= \overline{\Fr}_v$.  In examples  IIb${}_1$ and IIb${}_2$, 
with Frobenius classes $\Fr_v'$ and $\Fr_v''$ respectively, one has
 $\Fr'_{\sigma(v) }= {\Fr}''_v$.   Thus, for example, the first split prime
 $v=(1+2\gen)$ for $\Q(\sqrt{2})$ has conjugate $\sigma(v)=(3+2\gen)$,
 and in example Ia, one has $a_{\sigma(v)} = \alpha$, $d_{\sigma(v)} = \alpha^5$
 and  $\Fr_{\sigma(v)} = 10B$.  
 
 The agreement exhibited on Table~\ref{8matches} extends also to those $v$ with ${\mathbf N}(v)<100$ which do not have a place on the table.

\subsection*{Acknowledgements}  We are grateful to Laurent Berger for
 {the} suggestion of considering the Artin--Hasse exponential
in the context of a related question,  {and to David Savitt and Michael Schein
for discussions that confirmed the compatibility of Conjecture~\ref{conj:spaces}
with the results of \cite{gls}.
We thank Victor Abrashkin for calling our attention to the paper \cite{abrashkin},
and Toby Gee for informing us of work on Conjecture~\ref{conj:spaces} leading to its
proof in \cite{cegm}. We would also like to thank the referee for numerous suggestions that
improved the exposition of this paper.}

This research was partially supported by  EPSRC Grant EP/J002658/1 (LD),
Leverhulme Trust
RPG-2012-530,  EPSRC Grant  EP/L025302/1 and the Heilbronn Institute
for Mathematical Research (FD),  and
 Simons Foundation Collaboration Grant \#209472 (DPR).

\bibliographystyle{amsplain}
\bibliography{DDRSerreWild2.7.bib}

\providecommand{\bysame}{\leavevmode\hbox to3em{\hrulefill}\thinspace}
\providecommand{\MR}{\relax\ifhmode\unskip\space\fi MR }
\providecommand{\MRhref}[2]{%
  \href{http://www.ams.org/mathscinet-getitem?mr=#1}{#2}
}
\providecommand{\href}[2]{#2}
\begin{thebibliography}{10}

\bibitem{abrashkin}
Victor Abrashkin, \emph{Modular representations of the {G}alois group of a
  local field and a generalization of a conjecture of {S}hafarevich}, Izv.
  Akad. Nauk SSSR Ser. Mat. \textbf{53} (1989), no.~6, 1135--1182, 1337.
  \MR{1039960}

\bibitem{as}
Avner Ash and Warren Sinnott, \emph{An analogue of {S}erre's conjecture for
  {G}alois representations and {H}ecke eigenclasses in the mod {$p$} cohomology
  of {${\rm GL}(n,{\bf Z})$}}, Duke Math. J. \textbf{105} (2000), no.~1, 1--24.
  \MR{1788040}

\bibitem{breuil:appendix}
Christophe Breuil, \emph{Sur un probl\`eme de compatibilit\'e local-global
  modulo {$p$} pour {${\rm GL}_2$}}, J. Reine Angew. Math. \textbf{692} (2014),
  1--76. \MR{3274546}

\bibitem{bdj}
Kevin Buzzard, Fred Diamond, and Frazer Jarvis, \emph{On {S}erre's conjecture
  for mod {$\ell$} {G}alois representations over totally real fields}, Duke
  Math. J. \textbf{155} (2010), no.~1, 105--161. \MR{2730374}

\bibitem{cegm}
Frank Calegari, Matthew Emerton, Toby Gee, and Lambros Mavrides, \emph{Explicit
  {S}erre weights for two-dimensional {G}alois representations}, preprint,
  arxiv:1608.06059.

\bibitem{cd}
Seunghwan Chang and Fred Diamond, \emph{Extensions of rank one
  {$(\phi,\Gamma)$}-modules and crystalline representations}, Compos. Math.
  \textbf{147} (2011), no.~2, 375--427. \MR{2776609}

\bibitem{lassina:mathcomp}
Lassina Demb{\'e}l{\'e}, \emph{Explicit computations of {H}ilbert modular forms
  on {${\Bbb Q}(\sqrt{5})$}}, Experiment. Math. \textbf{14} (2005), no.~4,
  457--466. \MR{2193808}

\bibitem{lassina:expmath}
\bysame, \emph{Quaternionic {M}anin symbols, {B}randt matrices, and {H}ilbert
  modular forms}, Math. Comp. \textbf{76} (2007), no.~258, 1039--1057.
  \MR{2291849}

\bibitem{ddr}
Lassina Demb\'el\'e, Fred Diamond, and David~P. Roberts, \emph{On the
  computation and combinatorics of {S}erre weights for two-dimensional {G}alois
  representations}, in preparation.

\bibitem{doks}
Tim Dokchitser and Vladimir Dokchitser, \emph{Identifying {F}robenius elements
  in {G}alois groups}, Algebra Number Theory \textbf{7} (2013), no.~6,
  1325--1352. \MR{3107565}

\bibitem{fontaine}
Jean-Marc Fontaine, \emph{Il n'y a pas de vari\'et\'e ab\'elienne sur {${\bf
  Z}$}}, Invent. Math. \textbf{81} (1985), no.~3, 515--538. \MR{807070}

\bibitem{gee:type}
Toby Gee, \emph{Automorphic lifts of prescribed types}, Math. Ann. \textbf{350}
  (2011), no.~1, 107--144. \MR{2785764}

\bibitem{ghs}
Toby Gee, Florian Herzig, and David Savitt, \emph{{G}eneral {S}erre weight
  conjectures}, preprint, arxiv:1509.02527.

\bibitem{gk}
Toby Gee and Mark Kisin, \emph{The {B}reuil-{M}\'ezard conjecture for
  potentially {B}arsotti-{T}ate representations}, Forum Math. Pi \textbf{2}
  (2014), e1, 56. \MR{3292675}

\bibitem{gls}
Toby Gee, Tong Liu, and David Savitt, \emph{The {B}uzzard-{D}iamond-{J}arvis
  conjecture for unitary groups}, J. Amer. Math. Soc. \textbf{27} (2014),
  no.~2, 389--435. \MR{3164985}

\bibitem{jr}
John~W. Jones and David~P. Roberts, \emph{A database of number fields}, LMS J.
  Comput. Math. \textbf{17} (2014), no.~1, 595--618. \MR{3356048}

\bibitem{kw1}
Chandrashekhar Khare and Jean-Pierre Wintenberger, \emph{Serre's modularity
  conjecture. {I}}, Invent. Math. \textbf{178} (2009), no.~3, 485--504.
  \MR{2551763}

\bibitem{kw2}
\bysame, \emph{Serre's modularity conjecture. {II}}, Invent. Math. \textbf{178}
  (2009), no.~3, 505--586. \MR{2551764}

\bibitem{mavrides}
Lambros Mavrides, \emph{On wild ramification in reductions of two-dimensional
  crystalline {G}alois representations}, Ph.D. thesis, King's College London,
  2016.

\bibitem{newton}
James Newton, \emph{Serre weights and {S}himura curves}, Proc. Lond. Math. Soc.
  (3) \textbf{108} (2014), no.~6, 1471--1500. \MR{3218316}

\bibitem{robert}
Alain~M. Robert, \emph{A course in {$p$}-adic analysis}, Graduate Texts in
  Mathematics, vol. 198, Springer-Verlag, New York, 2000. \MR{1760253}

\bibitem{schein:ijm}
Michael~M. Schein, \emph{Weights in {S}erre's conjecture for {H}ilbert modular
  forms: the ramified case}, Israel J. Math. \textbf{166} (2008), 369--391.
  \MR{2430440}

\bibitem{serre_cl}
Jean-Pierre Serre, \emph{Local fields}, Graduate Texts in Mathematics, vol.~67,
  Springer-Verlag, New York-Berlin, 1979, Translated from the French by Marvin
  Jay Greenberg. \MR{554237}

\bibitem{serre}
\bysame, \emph{Sur les repr\'esentations modulaires de degr\'e {$2$} de {${\rm
  Gal}(\overline{\bf Q}/{\bf Q})$}}, Duke Math. J. \textbf{54} (1987), no.~1,
  179--230. \MR{885783}

\bibitem{rlt:inv}
Richard Taylor, \emph{On {G}alois representations associated to {H}ilbert
  modular forms}, Invent. Math. \textbf{98} (1989), no.~2, 265--280.
  \MR{1016264}

\end{thebibliography}

\end{document}